\documentclass[12pt,draftclsnofoot,onecolumn]{IEEEtran}
\usepackage{graphicx}
\usepackage{cite}
\usepackage{amsmath,amssymb,placeins}
\usepackage{algorithmic}
\usepackage{graphicx}
\usepackage{textcomp}
\usepackage[mathscr]{euscript}
\usepackage{lipsum}
\usepackage{cuted}
\usepackage{bm}
\usepackage{multicol}
\usepackage{array}
\usepackage{makecell}
\usepackage{diagbox}
\usepackage[colorlinks=true,linkcolor=gray,citecolor=black]{hyperref}
\usepackage[symbol]{footmisc}
\usepackage[utf8]{inputenc}
\usepackage{caption}
\usepackage{tabularx}

\usepackage[mathscr]{euscript}
\usepackage[usenames]{color}
\usepackage[table]{xcolor}
\usepackage{amsfonts}
\usepackage{latexsym}
\usepackage{subfigure}
\usepackage{tikz}
\usetikzlibrary{arrows.meta}
\usetikzlibrary{shapes.geometric, arrows}
\usepackage{caption}
\usepackage{floatrow}
\usepackage{multirow}
\usepackage{epstopdf}
\usepackage[english]{babel}

\newcommand{\ccr}{\color{red}}
\newcommand{\ccn}{\color{black}}
\newcommand{\ggr}{\color{gray}}

\newcommand{\mbb}{\mathbb}
\newcommand{\mbf}{\mathbf}
\newcommand{\om}{\zeta_q}
\newcommand{\mal}{\mathcal}

\newcolumntype{\arow}{\rightarrow}

\newcommand{\bx}{\mathbf{x}}
\newcommand{\be}{\mathbf{e}}
\newcommand{\bZ}{\mathbf{Z}}
\newcommand{\bbZ}{\mathbb{Z}}
\newcommand{\calV}{\mathcal{V}}
\DeclareMathOperator{\wt}{wt}
\usepackage{amsthm}
\newtheorem{theorem}{Theorem}
\newtheorem{lemma}{Lemma}
\newtheorem{corollary}{Corollary}

\newtheorem{remark}{Remark}

\newtheorem{example}{Example}

\def\BibTeX{{\rm B\kern-.05em{\sc i\kern-.025em b}\kern-.08em
		T\kern-.1667em\lower.7ex\hbox{E}\kern-.125emX}}

\def\@fnsymbol#1{\ensuremath{\ifcase#1\or *\or \dagger\or \ddagger\or
		\mathsection\or \mathparagraph\or \|\or **\or \dagger\dagger
		\or \ddagger\ddagger \else\@ctrerr\fi}}

\newcommand{\ba}{\mathbf{a}}
\newcommand{\bb}{\mathbf{b}}
\newcommand{\by}{\mathbf{y}}

{
\begin{document}
\title{A Further Investigation on Complete Complementary Codes from $q$-ary Functions}
\author{Palash Sarkar,~Chunlei Li,~Sudhan Majhi,~and ~Zilong Liu
\thanks{Palash Sarkar and Chunlei Li are with the Department of Informatics, Selmer Center, University of Bergen, Norway, e-mail: {\tt palash.sarkar@uib.no; chunlei.li@uib.no}.}% <-this % stops a space
\thanks{Sudhan Majhi is with the Department of Electrical Communication Engineering, Indian Institute of Science, Bangalore, India, e-mail:{\tt smajhi@iisc.ac.in}.}
\thanks{Zilong Liu is with the School of Computer Science and Electronics Engineering, University of Essex, Colchester CO4 3SQ, UK, e-mail:{\tt zilong.liu@essex.ac.uk}.}} 
\IEEEpeerreviewmaketitle	
\maketitle
\begin{abstract}
This research focuses on constructing $q$-ary functions for complete complementary codes (CCCs) with flexible parameters, such as length, set size, and alphabet. In recent years, significant effort has been dedicated to discovering new CCCs and their functional representations. However, most existing work has primarily identified sufficient conditions for $q$-ary functions related to $q$-ary CCCs. To the best of the authors' knowledge, this study is the first to establish both the necessary and sufficient conditions for $q$-ary functions, encompassing most existing CCCs constructions as special cases.
For $q$-ary CCCs with a length of $q^m$ and a set size of $q^{n+1}$, we begin by analyzing the necessary and sufficient conditions for $q$-ary functions defined over the domain $\mathbb{Z}_q^m$, where $q,~m \geq 2$ are integers, and $0 \leq n < m$. Additionally, we construct CCCs with lengths given by $L = \prod_{i=1}^k p_i^{m_i}$, set sizes given by $K = \prod_{i=1}^k p_i^{n_i+1}$, and an alphabet size of $\nu = \prod_{i=1}^k p_i$, where $p_1 < p_2 < \cdots < p_k$. Here, the parameters $m_1, m_2, \ldots, m_k$ are positive integers (each at least 2), and $n_1, n_2, \ldots, n_k$ are non-negative integers satisfying $n_i \leq m_i - 1$. The variable $k$ is a positive integer. To achieve these specific parameters, we examine the necessary and sufficient conditions for $\nu$-ary functions over the domain $\mathbf{Z}_{p_1}^{m_1} \times \cdots \times \mathbf{Z}_{p_k}^{m_k}$, which is a subset of $\mathbb{Z}_{\nu}^m$ and contains $\prod_{i=1}^k p_i^{m_i}$ vectors. In this context, $\mathbf{Z}_{p_i}^{m_i} = \{0, 1, \ldots, p_i - 1\}^{m_i}$, and $m$ is the sum of $m_1, m_2, \ldots, m_k$.
The $q$-ary and $\nu$-ary functions allow us to cover all possible length sequences. However, we find that the proposed $\nu$-ary functions are more suitable for generating CCCs with a length of $L = \prod_{i=1}^k p_i^{m_i}$, particularly when $m_i$ is coprime to $m_j$ for some $1 \leq i \neq j \leq k$. While the proposed $q$-ary functions can also produce CCCs of the same length $L$, the set size and alphabet size become as large as $L$, since in this case, the only choice for $q$ is $L$. In contrast, the proposed $\nu$-ary functions yield CCCs with a more flexible set size $K\leq L$ and an alphabet size of $\nu<L$. 
\end{abstract}
\section{Introduction}
\label{sec:intro}
\subsection{Background of CCCs}
The study of complementary sequences has been attracting increasing research attention in recent years due to their interesting
correlation properties as well as their wide applications in engineering.  Formally, a collection of sequences with identical
length and zero aperiodic auto-correlation sum (AACS) for all non-zero time shifts are called a complementary set (CS)
and each constituent sequence is called a complementary sequence \cite{chinchong,pater2000,liug}. In particular, a CS of size two
is called a Golay complementary pair (GCP), a widely known concept proposed by Marcel J. E Golay {in the late 1940s} \cite{Thesis_1949golay}.
{In this case,} the two constituent sequences are called Golay sequences.
Furthermore, by rearranging every CS {to be an ordered set}, a collection of CSs with zero aperiodic cross-correlation sum (ACCS)
for all time shifts is called mutually orthogonal Golay complementary sets (MOGCSs) \cite{rati}. It is noted that the set
size (denoted by $K$) of MOGCSs is upper bounded by the number of constituent sequences (denoted by $M$), i.e., $K\leq M$.
When $K=M$, the resultant MOGCSs are called a set of complete complementary codes (CCCs) \cite{hator}.

So far, complementary sequences and CCCs have been widely applied in numerous applications in coding, signal processing,
and wireless communications. For example, they {have been employed} for peak-to-mean envelope power ratio (PMEPR)
reduction in code-keying orthogonal frequency division multiplexing (OFDM) \cite{Davis1999},
signature sequences for interference-free multicarrier code-division multiple-access (MC-CDMA) over asynchronous wireless channels \cite{chen2007next,liumc}, training/sensing sequences in multiple-input
multiple-output (MIMO) communication/radar systems \cite{Wang2007,Pezeshiki2008,Tang2014},
information hiding in image/audio/video signals \cite{Thesis_2014Kojima},
kernel/seed sequences to construct Z-complementary code set (ZCCS) \cite{psktcom,pa_pbf} and quasi-complementary sequence set \cite{avikr_qccs, zhu_qccs,palash_qcss_tit_24} for MC-CDMA,
and zero-correlation zone (ZCZ) sequences \cite{dfan,ltsu,appus,Tang2010,Liu_ITW2014} in quasi-synchronous direct-sequence CDMA communications.

{There have been significant research attempts} on efficient constructions of complementary
sequences (including GCPs, MOGCSs, and CCCs) with flexible lengths. With interleaving and concatenation,
Golay first constructed binary GCPs of lengths $2MN$ from two shorter GCPs with lengths $M$ and $N$,
respectively  \cite{golay1961}. Budi\v{s}in developed recursive algorithms for more polyphase and multi-level
complementary sequences through a series of sequence operations (e.g., phase rotations, shifting) associated
with permutation vectors \cite{Budisin90a,Budisin90b}. In 1999, Davis and Jedwab proposed a systematic
construction of $2^h$-ary ($h\in\mathbb{Z}^+$) GCPs with power-of-two lengths using $2$nd-degree generalized
Boolean functions (GBFs) \cite{Davis1999}. Such Golay sequences are called Golay-Davis-Jedwab (GDJ) {complementary} sequences in this paper.
Paterson extended the idea of $2^h$-ary GCPs to $q$-ary (even $q$) CSs with power-of-two constituent sequences {\color{red}and} power-of-two
lengths \cite{pater2000}. His work was further generalized in 2007 by Schmidt by moving to the higher-order GBFs  \cite{Schmid2007}.
Recently, a construction is reported in \cite{palcs} using higher-order GBFs for complementary sequences with
lower PMEPR upper bound. For complementary sequences with non-power-of-two lengths, Chen contributed with the aid of certain truncated
Reed-Muller codes \cite{chentit}.

The constructions of MOGCSs especially CCCs are more stringent owing to the additional requirement of zero ACCS  \cite{rati,liumc,chencommlett}. In 2020, Wu-Chen-Liu introduced a construction of MOGCSs with non-power-of-two lengths based on GBFs in \cite{swuc} where the set size is exactly half of the number of constituent sequences in a CS, i.e., $K=M/2$. Later several attempts have been made to construct non-power-of-two lengths MOGCSs with parametric restriction $K=M/2$, such as in \cite{nwmocs1,xiao2023new}. However, in a recent work \cite{pku1}, the authors successfully constructed binary CCCs of lengths in the form of $5\cdot 2^m$ and $13\cdot 2^m$, and the set size in the form of power-of-two. In \cite{shen2023}, $q$-ary extended Boolean functions (EBFs) based construction of CCCs of length $q^m$ has been reported, where $q\geq 2$ is any integer. In \cite{sarkar2021multivariable}, we proposed construction of CCCs, of length $\prod_{i=1}^kp_i^{m_i}$ and set size $p_i^{n_i+1}$, over an alphabet of size $q=\prod_{i=1}^kp_i$, where the prime numbers $p_i$ are distinct, and $0\leq n_i<m_i-1$. The construction reported in \cite{sarkar2021multivariable} is based on $q$-ary functions which we describe in the below subsection. Another powerful tool is paraunitary (PU) generator with which a set of CCCs can be efficiently obtained by sending a pipeline of shifted impulses to multistage filterbanks whose coefficients are extracted from a group of PU matrices \cite{Budisin_QAM,Wang_SETA2016,sdas,Sdas_lett,shibu2}. Very recently, the inherent connections between the GBF generator and the PU generator in CCCs construction have been uncovered in \cite{wang2020new} and in connection with this work, Wang-Gong further proposed CCCs with the help of Buston-type Hadamard matrices and two-level auto-correlation functions, in \cite{wangong} where they also established a theoretical connection between the constructed CCCs and their functional representations. 
In this paper, we shall propose functions with two different kinds of domains. To avoid confliction, instead of using the terms GBF or EBF, we call a function as $q$-ary function if its co-domain is $\mbb Z_q$.  % when some good PU matrix kernels can be found.
\subsection{Contributions and Approach}
In most of the existing construction of CCCs of lengths in the form $q^m$ and set sizes in the form $q^{n+1}$, where $q\geq 2$ is any integer and $0\leq n<m$, each sequence can be defined by a $q$-ary function $f:\mbb Z_q^m\rightarrow \mbb Z_q$, where the restriction of $f$ at $\mbf x_J=\mbf c\in \mbb Z_q^n$ can be expressed as 
$f\arrowvert_{\mathbf{x}_J=\mathbf{c}}=\sum_{i=1}^{m-n-1} h_i(x_{\pi_c(i)})h'_i(x_{\pi_c(i+1)})+\sum_{j=1}^{m-n} g_j(x_{\pi_c(j)})+c$. Based on the existing works, many different forms for the functions $h_i,~h_i',~g_j: \mbb Z_q\rightarrow\mbb Z_q$ have already been reported. Besides, we define  $J=\{j_1,j_2,\hdots,j_n\}\subset \{1,2,\hdots,m\}$, $\mbf x_J=(x_{j_1},x_{j_2},\hdots, x_{j_n})$, $\mbf c=(c_1,c_2,\hdots,c_k)~\in\mbb Z_q^k$, $c=\sum_{i=1}^k c_i q^{k-i} $, and $\pi_c:\{1,2,\hdots,m-n\}\rightarrow \{1,2,\hdots,m\}\setminus J $ is a one-to-one mapping. The construction of CCCs in \cite{rati}, sequences can be expressed by $f$, where  $h_i(x_{\pi_c(i)})=x_{\pi_c(i)}$, $h'_i(x_{\pi_c(i+1)})=x_{\pi_c(i+1)}$, and $q=2$. This form of $f$ in \cite{rati}, produces binary CCCs of lengths $2^m$ and set size $2^{n+1}$. The smallest alphabet size that can be obtained from \cite{rati} is $2$, however, this work also produces CCCs over any even size alphabet. In 2023 \cite{wangong}, Wang-Gong proposed a generic functional representation for CCCs of lengths $p^m$ and set size $p$, where $p$ is a prime number. The function $f:\mbb Z_p^m\rightarrow \mbb Z_p$ in \cite{wangong} can be obtained by substituting $n=0$, $J=\emptyset$, and $q=p$, in $f\arrowvert_{\mathbf{x}_J=\mathbf{c}}$, where $h_i,~h_i'$ are any functions that permute $\mbb Z_p$. Later that same year, in pursuit of CCCs with a length of $q^m$ and a set size of $q^{n+1}$, \cite{shen2023} introduced certain specific $q$-ary functions, where $h_i(x_{\pi_c(i)})=a_i x_{\pi_c(i)}$, $h'_i(x_{\pi_c(i+1)})=x_{\pi_c(i+1)}$, $g_j(x_{\pi_c(j)})=\sum_{i=1}^{q-1} c_{i,j}x_{\pi_c(j)}^i+c_j$, $a_i$ is an integer co-prime to $q$, and $c_{i,j},~c_j\in \mbb Z_q$. In much of the existing research on functional representations of CCCs, the functions $h_i$ and $h_i'$ are typically selected as specific functions that notably permute $\mathbb{Z}_q$. 

The aforementioned matters motivate us to study the following question: \textit{given a  CCCs of length $q^m$ and set size $q^{n+1}$, where each sequence can be explained by $f$ for which the restriction $f\arrowvert_{\mathbf{x}_J=\mathbf{c}}$ is given by the aforementioned expression, what will be the properties of the functions $h_i$, $h_i'$ and $g_j$?} To analyze this, we dive into the properties of CCCs i.e., the out-of-phase aperiodic auto/cross-correlation functions of the codes in a CCCs are zero. We have shown that the properties of CCCs imply that the function $h_i$ and $h_i'$ must permute $\mbb Z_q$, and $g_j$ can be any function on $\mbb Z_q$. To derive the aforementioned properties of $h_i$ and $h_i'$, we apply the zero correlation property between two CSs at a time shift $\tau=q^{m-n}-rq^{m-n-i}$, where $1\leq r<q$. By choosing  $r=1$, we have shown that the zero auto/cross-correlation property between two CSs at $\tau=q^{m-n}-q^{m-n-i}$ implies $h_i$ permutes $\mbb Z_q$. Then by going through the rest of the values of $r$, we reach to the condition that $h_i'$ also permute $\mbb Z_q$. Unlike the existing constructions, this approach provides a complete classification, i.e., both of necessary and sufficient conditions for the afore-mentioned $q$-ary functions $f$, and provides a large class of functions that includes not only most of the existing works as a special case but also produces many CCCs that are unlikely be covered in existing works. 
  
The aforementioned functions can generate CCCs with length $q^m$ for any positive integer $q\geq 2$. 
 However, to obtain CCCs of length $L$ that is not a power of $q$, 
it can be observed that one needs to choose the alphabet size as $L$, and so the set size also becomes as large as $L$. CCCs with such parameters can be derived from the DFT matrix as well \cite{jinla}.
Such CCCs may not be feasible for use scenarios where low alphabet size and small set size are required. For example, in machine-type communications, we need to maintain a small alphabet size, and for OFDM application to deal with high PMEPR, one needs sequences coming from a small size CS. To overcome this limitation, we introduce another type of $q$-ary functions $f$ on the domain  $\mbf{Z}_{p_1}^{m_1}\times\mbf{Z}_{p_2}^{m_2}\times \cdots\times \mbf{Z}_{p_k}^{m_k}\subset \mbb Z_q^m$ and co-domain $\mathbb{Z}_{q}$, where $q=\prod_{i=1}^k p_i$, $m=\sum_{i=1}^k m_i$, and $\mbf Z_{p_i}=\{0,1,\hdots,p_1-1\}\subset \mbb Z_q$. On this setting of domain and co-domain, by analyzing both necessary and sufficient conditions, we have proposed a complete class of functions for CCCs with the aforementioned parameters $K$, $L$, and $q$, in \textbf{Corollary \ref{corr19824}}. To keep the presentation simple, let us explain our contribution for $k=2$ and $n_1=n_2=0$. In this case, the proposed functions appear in the following form $f(x_1,x_2,\hdots,x_m)=(\frac{q}{p_1}\sum_{i=1}^{m_1-1} f_i(x_{\pi(i)})f_i'(x_{\pi(i+1)})+\sum_{\alpha=1}^{m_1}g_{\alpha}(x_\alpha))+(\frac{q}{p_2}\sum_{j=1}^{m_2-1} h_j(x_{\pi'(m_1+j)}) h_j'(x_{\pi'(m_1+j+1)})+\sum_{\beta=1}^{m_2}g'_{\beta} ({x_{\beta+m_1}}))+\gamma f_0(x_{\pi(m_1)})h_0(x_{\pi'(1)})$, where $f_i,~f_i'$, $h_j,~h_j',~f_0,~g_{\alpha},~h_0,~g'_{\beta}$ are functions from $\mbb Z_q$ to $\mbb Z_q$, and $\pi:\{1,2,\hdots,m_1\}\rightarrow\{1,2,\hdots,m_1\}$ and $\pi':\{m_1+1,m_1+2,\hdots,m\}\rightarrow \{m_1+1,m_1+2,\hdots,m\}$ are one-to-one mappings. We have shown that for any choice of the functions $f_0,h_0,g_\beta'$, and one-to-one mappings $\pi,\pi'$, a set of codes constructed from $f$ will form CCCs iff $f_i$ and $f_i'$ permute $\mbf Z_{p_1}$ under modulo $p_1$ and $h_j,h_j'$ permute $\mbf Z_{p_2}$ under modulo $p_2$. To obtain this property of $f_i$ and $f_i'$, we made more detailed investigation of correlation property at different time shifts. Specifically, we use the zero correlation property between two CSs at the time shift $\tau=p_1^{m_1}-np_1^{m_1-i}$. By choosing $n=1$, we have shown that the zero correlation property at the time shift $\tau=p_1^{m_1}-p_1^{m_1-i}$ implies $f_i$ permutes $\mbf Z_{p_1}$ under modulo $p_1$ operation. Then by going through the rest of the values of $n$, it has been shown that $f_i'$ also permute $\mbf Z_{p_1}$ under modulo $p_1$. We apply a similar approach to show that $h_j,h_j'$ permute $\mbf Z_{p_2}$ under modulo $p_2$. In \cite{sarkar2021multivariable}, the proposed $q$-ary function can be obtained by choosing $f_i(x_{\pi(i)})=x_{\pi(i)}$, $f_i'(x_{\pi(i+1)})=x_{\pi(i+1)}$, $ h_j(x_{\pi'(m_1+j)})=x_{\pi'(m_1+j)}$, $h_j'(x_{\pi'(m_1+j+1)})=x_{\pi'(m_1+j+1)}$, $g_{\alpha}(x_\alpha)=c_\alpha x_\alpha+c$, and $g'_{\beta} ({x_{\beta+m_1}})=c_{\beta+m_1}x_{\beta+m_1}+c'$, where $c_\alpha, c_{\beta+m_1}, c,c'\in \mbb Z_q$.
As every integer greater than $1$ can be factored
uniquely as the product {of powers of} prime numbers (by unique factorization theorem), it can be observed that the proposed $q$-ary functions on the domain $\mbf{Z}_{p_1}^{m_1}\times\mbf{Z}_{p_2}^{m_2}\times \cdots\times \mbf{Z}_{p_k}^{m_k}$ is also capable of producing CCCs of arbitrary length. 

The structure of this paper is organized as follows: Section \ref{premil20924} introduces the essential mathematical tools used in this study. In Section \ref{sectionIII200924}, we examine the necessary and sufficient conditions for $q$-ary functions on the domain $\mbb Z_q^m$ in relation to CCCs. Section \ref{sectn4} extends this analysis to $q$-ary functions on the domain $\mbf{Z}_{p_1}^{m_1}\times\mbf{Z}_{p_2}^{m_2}\times \cdots\times \mbf{Z}_{p_k}^{m_k}$, also in connection with CCCs. Finally, Section \ref{conclu200924} concludes the paper.
\section{Preliminaries}\label{premil20924}
 In this paper, the following notations are used otherwise stated. 
\begin{itemize}
	\item For a finite set $S$, $|S|$ denotes the cardinality of $S$.
	\item For two finite sets $S$ and $S'$, $S\times S'$ denotes the Cartesian product between $S$ and $S'$.
	\item $S^n=S\times \cdots\times S$ indicates the Cartesian product of $n$ identical sets $S$.
    \item $\mathbb{Z}$ denotes the set of all integers.
    \item $\mathbb{Z}_q$ denotes the set of all integers modulo $q$.
	\item $\mathbf{Z}_t=\{0,1,\hdots,t-1\}\subset \mathbb{Z}_q $ for a positive integer $t|q$.
    \item $L$ is a positive integer having factorization as 
    $L=\prod_{i=1}^kp_i^{m_i}$ with $k\geq 2$, $m_i\geq 2$, where  $p_1<\cdots <p_k$.
    \item $K=\prod_{i=1}^k p_i^{n_i+1}$, where $0\leq n_i\leq m_i-1$ is a non-negative integer for $1\leq i\leq k$.
    \item $q = \prod_{i=1}^k p_i$ and $m=\sum_{i=1}^km_i$.
    \item $\mal{V}_L=\mbf{Z}_{p_1}^{m_1}\times\mbf{Z}_{p_2}^{m_2}\times\cdots\times \mbf{Z}_{p_k}^{m_k}\subset \mathbb{Z}_q^m$.
    \item $\xi_q=\exp(2\pi\sqrt{-1}/q)$ is a primitive $q$-th root of unity.
    \item $\mathcal{A}_q = \{\xi_q^i\,:\, i = 0,1,\dots, q-1\}$.
	%\item $q$ is a positive integer such that $p_i|q$, where $i=1,2,\hdots,k$.
    \item For two vectors $\ba, \bb$ of length $N$, their inner product is denoted by $\ba\cdot \bb$.
\end{itemize} \ccn
\subsection{Aperiodic Auto- and Cross-Correlation}
%Assume that $\tau$ is an integer satisfying $0\leq |\tau|<L$.
For two complex-valued sequences $\mathbf{a}=(a_0,a_1,\hdots,a_{L-1})$ and $\mathbf{b}=(b_0,b_1,\hdots,b_{L-1})$ of length $L$, we define the aperiodic cross-correlation function (ACCF) at a shift $\tau$, with  $0\leq |\tau|<L$ as
\begin{equation}\label{accf}
\Theta (\mathbf{a},\mathbf{b})(\tau)=
\begin{cases}
\sum_{t=0}^{L-\tau-1}a_t b^*_{t+\tau}, & 0\leq \tau<L,\\
\sum_{t=0}^{L+\tau-1}a_{t-\tau} b^*_t,&-L<\tau<0.
\end{cases}
\end{equation}
For $\mathbf{a}=\mathbf{b}$, the ACCF defined in (\ref{accf}) reduces to the aperiodic auto-correlation function (AACF) of 
$\mathbf{a}$, which will be denoted as $\Theta (\mathbf{a})(\tau)$ for short. 

Let $\mathcal{C}$ be a $K\times M$ matrix, of which each entry is a sequence of length $L$. We denote by $C_k$ the $k$-th row of sequences in $\mal C$, given by 
$$\mathbf{C}_k=\begin{bmatrix}
\mathbf{c}_{k}^1&
\mathbf{c}_{k}^2&
\cdots&
\mathbf{c}_{k}^M
\end{bmatrix},$$
with $\mathbf{c}_{k}^m$ being a complex-valued sequence of length $L$, where $1\leq m\leq M$.
The ACCF (sum) between $\mathbf{C}_{k_1}$
and $\mathbf{C}_{k_2}$ for $1\leq k_1, k_2 \leq K$ is defined as 
\begin{equation}\label{accfk}
\Theta (\mathbf{C}_{k_1},\mathbf{C}_{k_2})(\tau)=\sum_{m=1}^{M}\Theta   (\mathbf{c}_{k_1}^m,\mathbf{c}_{k_2}^m)(\tau),
\end{equation}
For $k_1=k_2=k$, the ACCF in (\ref{accfk}) reduces to the AACF of $\mathbf{C}_{k}$ and we denote it by $\Theta (\mathbf{C}_{k})(\tau)$. The set $\mal{C}$ is said to be a mutually orthogonal Golay complementary set (MOGCS) if
it satisfies the following properties:
\begin{equation}\label{accfc}
\begin{split}
\Theta (\mathbf{C}_{k_1},\mathbf{C}_{k_2})(\tau)=
\begin{cases}
ML,& \tau=0,~k_1=k_2,\\
0,& 0<|\tau|<L, ~k_1=k_2,\\
0,&|\tau|<L, ~k_1\neq k_2.
\end{cases} 
\end{split}
\end{equation}
The above property can be represented by polynomials more compactly. 
Given a sequence $\ba=(a_0,a_1,\dots, a_{L-1})$, its generating polynomial is given by $A(z)=\sum_{i=0}^{L-1}a_iz^i$.  Note that for any two sequences $\ba$ and  $\bb$, their generating polynomials satisfy
$$A(z)\overline{B}(z^{-1})= \sum_{|\tau| <L} \Theta (\mathbf{a},\mathbf{b})(\tau) z^{\tau}, $$ where $\overline{B}(z)$ is the conjugate polynomial of $B(z)$ given by $\overline{B}(z) = \sum_{i=0}^{L-1} b^*_i z^i$. 
A code set $\mathcal{C}=\{\mbf{C}_1,\dots, \mathbf{C}_K\}$ can be represented as a $K\times M$ matrix $\mathbf{C}(z)$
whose $(k,m)$-th entry is the associated polynomial of the sequence $\mathbf{c}_k^m$. 
Then the code set $\mathcal{C}$ is an MOGCS 
if 
\[
\mathbf{C}(z)\cdot \mathbf{C}^{\dagger}(z^{-1}) = ML\cdot I_K,
\]
where $\mathbf{C}^{\dagger}(z^{-1})$ is the conjugate transpose matrix of $\mathbf{C}(z^{-1})$ by and $I_K$ is the identity matrix of order $K$. Since $I_K$ has rank $K$,
the matrix $\mbf{C}(z)$ has rank at least $K$, implying that $K\leq M$ for any MOGCS \cite{schweitzer1971, rati}. 
In particular, when $K=M$, a MOGCS is called a set of complete complementary codes (CCCs).
In this case, we denote the code set $\mathbf{C}$ as $(M,L)$-CCCs, and any of the codes in $\mal{C}$ is known as a CS or complementary code (CC) \cite{sarkar2021multivariable}. 
\begin{lemma}\label{nemma}
Let $q$ be a positive integer and $\phi$ be a function from $\mbb Z_q$ to itself. Then 
\begin{equation}
\begin{split}
\sum_{x\in \mbb Z_q} \xi_q^{r\cdot \phi(x)}=0  \text{ for any } r\in \mbb Z_q\setminus\{0\}~\text{if and only if}~\phi~\text{is a permutation of}~\mbb Z_q.
\end{split}
\end{equation}
\end{lemma}
\begin{IEEEproof} Denote by $\chi_r(x) = \xi_q^{r\cdot x}$, the function from $\mbb Z_q$ to $\mathcal{A}_q$. It is easily seen that $\chi_r = \xi_q^r$ for $r=0,1,\dots, q-1$
are (additive) characters of $\mbb Z_q$. Recall from \cite[Sec. 5.1]{Lidl_Niederreiter_1996} that 
\begin{equation}\label{eq_char}
\sum_{x\in \mbb Z_q} \chi_r(x) = 0, \,\, \forall \, r\in \mbb Z_q\setminus\{0\},
\end{equation}
which implies 
\begin{equation}\label{eq_char_1}
\sum_{x\in \mbb Z_q} \chi_{r_1}(x) \overline{\chi_{r_2}(x)} = \sum_{x\in \mbb Z_q} \chi_{r_1-r_2}(x) = 
\begin{cases}
    q, & \text{if } r_1 = r_2, \\
    0, & \text{if } r_1\neq r_2,    
\end{cases}
\end{equation}
and 
\begin{equation}\label{eq_char_2}
\sum_{r\in \mbb Z_q} \chi_{r}(x) \overline{\chi_{r}(y)} = 
\sum_{r\in \mbb Z_q} \chi_{r}(x-y) = 
\begin{cases}
    q, & \text{if } x = y, \\
    0, & \text{if } x\neq y.    
\end{cases}
\end{equation}
When $\phi$ permutes $\mbb Z_q$, it follows from \eqref{eq_char} that
$$
\sum_{x\in \mbb Z_q} \xi_q^{r\cdot \phi(x)} = \sum_{x\in \mbb Z_q} \xi_q^{r\cdot x}, \,\, \forall \, r\in \mbb Z_q\setminus\{0\}.
$$
Conversely, for any $a\in \mbb Z_q$, the number $N=\#\{c\in \mbb Z_q\,:\, \phi(x) = a\}$
is given by
\[
\begin{split}
    N = \frac{1}{q}\sum_{c\in \mbb Z_q}\sum_{r\in \mbb Z_q} \chi_r(\phi(c)) \overline{\chi_r(a)} = \frac{1}{q}\sum_{r\in \mbb Z_q}\left(\sum_{c\in \mbb Z_q} \chi_r(\phi(c))\right) \overline{\chi_r(a)} = \frac{1}{q}\sum_{c\in \mbb Z_q} \chi_0(\phi(c))\chi_0(a) = 1,
\end{split}
\]which implies that $\phi$ permutes $\mbb Z_q$.
\end{IEEEproof}

\section{Construction of CCCs Using $q$-ary Functions }\label{sectionIII200924}
In this section, we are interested in finding a relationship between $q$-ary functions and CCCs for any positive integer values of $q\geq 2$. For a $q$-ary function $f:\mbb{Z}_q^m\rightarrow\mbb{Z}_q$ of $m$ variables $x_1,x_2,\hdots,x_m$, a complex-valued sequence over $\mathcal{ A}_q=\{\xi_q^i:i=0,1,\hdots,q-1\}$, can be defined as follows:
$$\psi(f)=(\zeta_q^{f_0},\zeta_q^{f_1},\hdots,\zeta_q^{f_{q^m-1}}),$$
where $f_x=f(\mbf x)$, $\mbf x=(x_1,x_2,\hdots,x_m)$ is the $q$-ary vector representation of $x$, i.e., $x=\sum_{j=1}^m x_j q^{j-1}$. It is to be noted that only when $q$ is a prime power,  one can obtain a unique algebraic expression of the function $f$ from the sequence $(f_0,f_1,\dots, f_{q^m-1})$. For $i=1,2,\hdots,m-1$ and $j=1,2,\hdots,m$, let us assume $h_i,~h_i'$ and $g_j$ are univariate functions from $\mbb Z_q$ to $\mbb Z_q$. Then for a permutation $\pi:\{1,2,\hdots,m\}\rightarrow \{1,2,\hdots,m\}$, define a $q$-ary function $f$ as follows:
\begin{equation}\label{pfunc_12724}
	\begin{split}
	f(x_1,x_2,\hdots,x_m)=\sum_{i=1}^{m-1} h_i(x_{\pi(i)}) h_i'(x_{\pi(i+1)})+\sum_{j=1}^m g_j(x_j).
	\end{split}
\end{equation} 
For $t\in\mbb Z_q$, define 
\begin{equation}\nonumber
	\begin{split}
     C_t=\left\{f+dx_{\pi(1)}+tx_{\pi(m)} : d\in \mbb Z_q  \right\},
	\end{split}
\end{equation}
and
\begin{equation}\label{psi12724}
\psi(C_t)=\left\{ \psi\left(f+dx_{\pi(1)}+tx_{\pi(m)} \right):d\in\mbb Z_q \right\}.
\end{equation}
In the literature, several works have been reported on the construction of CCCs using $q$-ary functions \cite{rati,shen2023} where the main contribution is to find classes of $q$-ary functions that generate CCCs with flexible parameters. As an example, consider $h_i(x_{\pi(i)})=a_ix_{\pi(i)}$, $h_i'(x_{\pi(i+1)})=x_{\pi(i+1)}$, and $g_j(x_j)=\sum_{i=1}^{q-1}c_{i,j}x^i+c_j$ in (\ref{pfunc_12724}), where $a_i$ is an integer co-prime to $q$. Then substituting $f$ in (\ref{psi12724}), we obtain the $(q,q^m)$-CCCs in \cite{shen2023}. 

Besides CCCs, setting $q=2$, $h_i(x_{\pi(i)})=x_{\pi(i)}$ and $h_i(x_{\pi(i+1)})=x_{\pi(i+1)}$ in (\ref{pfunc_12724}), for $t=0$, the code defined in (\ref{psi12724}) forms binary GCP in \cite{Davis1999}. With the same setting as  described for \cite{Davis1999}, varying $t\in \mbb Z_2$ in (\ref{psi12724}), one can obatin binary $(2,2^m)$-CCCs from \cite{rati}. 

Interestingly, it can be observed that in the proposed function $f$ of \cite{rati,shen2023}, the functions $h_i$ and $h_i'$ are chosen in such a way that they permute $\mbb Z_q$. To the best of the author's knowledge, in most of the constructions for $(q,q^m)$-CCCs, the constituent sequences in each CSs can be presented by the function in (\ref{pfunc_12724}), provided that $h_i$ and $h_i'$ are some permutations of $\mbb Z_q$. In \cite{wangong}, it has been shown by Wang and Gong that the set of codes $\mal C=\{\psi(C_t): t\in \mbb Z_q \}$ forms a $(q,q^m)$-CCCs if the functions $h_i$ and $h_i'$ permute $\mbb Z_q$. It motivates us to further study $h_i$ and $h_i'$, mainly the required necessary condition for $h_i$ and $h_i'$, specially when $q\geq 2$ is any integer, so that the set of codes $\mal C=\{\psi(C_t): t\in \mbb Z_q \}$ forms $(q,q^m)$-CCCs.  
\ccn
\begin{theorem}\label{th824}
For any choice of the functions $g_j:\mbb Z_q\rightarrow \mbb Z_q$, $j=1,2,\hdots,m$, and permutation $\pi$, the set of complex-valued codes $\mathcal{C}=\{\psi(C_t):t\in\mbb Z_q \}$, where $\psi(C_t)$ is defined in (\ref{psi12724}), forms $(q,q^{m})$-CCCs iff  the functions $h_i$ and $h_i'$ permute $\mbb Z_q$, where $i=1,2,\hdots,m-1$.
\end{theorem}
\begin{IEEEproof}
Assume $L=q^m$, $x,y\in\mbb Z_L$, such that $y=x+\tau$ for some $\tau\geq 0$. Then for any two integers $t$ and $t'$ in $\mbb Z_q$, the ACCF between $\psi(C_t)$ and $\psi(C_{t'})$ at $\tau$ is given by 
\begin{equation}\label{key1_12724}
	\begin{split}
	\Theta(\psi(C_t),\psi(C_{t'}))(\tau)&=\sum_{d\in\mbb Z_q}\sum_{x=0}^{L-\tau-1} \xi_q^{F_{x,t}-F_{y,t'}+d(x_{\pi(1)}-y_{\pi(1)})}\\
	&=\sum_{x=0}^{L-\tau-1} \xi_q^{F_{x,t}-F_{y,t'}} \left(\sum_{d\in\mbb Z_q} \xi_q^{d(x_{\pi(1)}-y_{\pi(1)})}\right),
	\end{split}
\end{equation}
where $F_{x,t}=f(\bx)+tx_{\pi(m)}$, $F_{y,t'}=f(\by)+t'y_{\pi(m)}$. Now
\begin{equation}\label{key2_12724}
	\begin{split}
\sum_{d\in\mbb Z_q} \xi_q^{d(x_{\pi(1)}-y_{\pi(1)})}=
\begin{cases}
q, & x_{\pi(1)}=y_{\pi(1)},\\
0,& x_{\pi(1)}\neq y_{\pi(1)}.
\end{cases}
	\end{split}
\end{equation}
Then applying (\ref{key2_12724}) in (\ref{key1_12724}), we have
\begin{equation}\nonumber
	\begin{split}
\Theta(\psi(C_t),\psi(C_{t'}))(\tau)=q \sum_{\substack{x=0\\ x_{\pi(1)}=y_{\pi(1)}}}^{L-\tau-1} \xi_q^{F_{x,t}-F_{y,t'}}.
	\end{split}
\end{equation}
Let us assume $1<v\leq m$ be the smallest positive integer such that $x_{\pi(v)}\neq y_{\pi(v)}$. Let us now define two integers $x^u$ and $y^u$ in such a way that their vector representations differ from the vector representations of $x$ and $y$, respectively, only at the position $\pi(v-1)$, as follows:
$$(x_1,x_2,\hdots,(x_{\pi(v-1)}+u)\!\!\!\!\!\mod q,\hdots,x_m)~\text{and}~(y_1,y_2,\hdots, (y_{\pi(v-1)}+u)\!\!\!\!\!\mod q,\hdots,y_m).$$
Then
\begin{equation}\nonumber
	\begin{split}
F_{x^u,t}-F_{x,t}=&h_{v-2}(x_{\pi(v-2)})h_{v-2}'((x_{\pi(v-1)}+u)\!\!\!\!\!\mod q)+h_{v-1}((x_{\pi(v-1)}+u)\!\!\!\!\!\mod q)h_{v-1}'(x_{\pi(v)})\\
&-h_{v-2}(x_{\pi(v-2)})h_{v-2}'(x_{\pi(v-1)})-h_{v-1}(x_{\pi(v-1)})h_{v-1}'(x_{\pi(v)})\\
&+g_{\pi(v-1)}((x_{\pi(v-1)}+u)\!\!\!\!\!\mod q)-g_{\pi(v-1)}(x_{\pi(v-1)}),
	\end{split}
\end{equation}
and
\begin{equation}\nonumber
\begin{split}
F_{y^u,t'}-F_{y,t'}=&h_{v-2}(x_{\pi(v-2)})h_{v-2}'((x_{\pi(v-1)}+u)\!\!\!\!\!\mod q)+h_{v-1}((x_{\pi(v-1)}+u)\!\!\!\!\!\mod q)h_{v-1}'(y_{\pi(v)})\\
&-h_{v-2}(x_{\pi(v-2)})h_{v-2}'(x_{\pi(v-1)})-h_{v-1}(x_{\pi(v-1)})h_{v-1}'(y_{\pi(v)})\\
&+g_{\pi(v-1)}((x_{\pi(v-1)}+u)\!\!\!\!\!\mod q)-g_{\pi(v-1)}(x_{\pi(v-1)}).
\end{split}
\end{equation}
Then
\begin{equation}\nonumber
	\begin{split}
F&_{x^u,t}-F_{y^u,t'}-(F_{x,t}-F_{y,t'})\\&=F_{x^u,t}-F_{x,t}-(F_{y^u,t'}-F_{y,t'})\\
&=(h_{v-1}'(x_{\pi(v)})-h_{v-1}'(y_{\pi(v)}))(h_{v-1}((x_{\pi(v-1)}+u)\!\!\!\!\!\mod q)-h_{v-1}(x_{\pi(v-1)}))\\
&=rp(u),
	\end{split}
\end{equation}
where $r=h_{v-1}'(x_{\pi(v)})-h_{v-1}'(y_{\pi(v)})$ and $p(u)=h_{v-1}((x_{\pi(v-1)}+u)\!\!\!\mod q)-h_{v-1}(x_{\pi(v-1)})$.
Then
\begin{equation}\label{newcond_12724}
	\begin{split}
	\sum_{u=1}^{q-1} \xi_q^{F_{x^u,t}-F_{y^u,t'}-(F_{x,t}-F_{y,t'})}=-1+\sum_{u=0}^{q-1} \xi_q^{rp(u)}.
	\end{split}
\end{equation}

\noindent \textit{\textbf{\textbf{Sufficiency:}}}
We assume $h_i$ and $h_i'$ permute the set $\mbb Z_q$, where $i=1,2,\hdots,m-1$. Then 
 $r\neq 0$ and $p(u)$ also permutes $\mbb Z_q$, and
\begin{equation}\label{newcond1_12724}
\begin{split}
\xi_q^{F_{x,t}-F_{y,t'}}+\sum_{u=1}^{q-1} \xi_q^{F_{x^u,t}-F_{y^u,t'}}=0.
\end{split}
\end{equation}
Now from (\ref{newcond1_12724}), and combining all the above results of this proof, it is clear that 
\begin{equation}\label{newcond2_12724}
	\Theta(\psi(C_t),\psi(C_{t'}))(\tau)=0~\text{if}~h_{v-1},~ h_{v-1}'~\text{permute}~\mbb Z_q,
\end{equation}
where $v=2,3,\hdots,m$. 
Applying a similar technique as $\tau>0$, we can reach the same result as in (\ref{newcond2_12724}) for $\tau<0$. 
Now for $\tau=0$, we have
\begin{equation}\nonumber
	\begin{split}
	\Theta(\psi(C_t),\psi(C_{t'}))(0)=q\sum_{x=0}^{L-1} \xi_q^{(t-t')x_{\pi(m)}}=
	\begin{cases}
	q^{m+1},& t=t',\\
	0,&\text{otherwise}.
	\end{cases}
	\end{split}
\end{equation}
Hence the set of codes $\mathcal{C}$ forms $(q,q^m)$-CCCs if $h_i$ and $h_i'$ permute $\mbb Z_q$, where $i=1,2,\hdots,m-1$.

\noindent \textit{\textbf{Necessity:}}
Let us assume that, for any choice of the functions $g_j:\mbb Z_q\rightarrow \mbb Z_q$, $j=1,2,\hdots,m$, and one-to-one mapping $\pi$ in (\ref{pfunc_12724}), the set of complex-valued codes $\mathcal{C}=\{\psi(C_t):t\in\mbb Z_q \}$ forms $(q,q^{m})$-CCCs. Using this property,  our goal is to determine the bijectivity of the functions $h_i$ and $h_i'$, where $i=1,2,\hdots,m-1$. Let us assume $\tau=q^m-nq^{m-i}$, and $\pi(i')=i'$, where $1\leq n\leq q-1$, and $i'=1,2,\hdots,m$. As $\mathcal{C}$ forms CCCs, 
then 
\begin{equation}\nonumber
	\begin{split}
	\Theta(\psi(C_t),\psi(C_{t'}))(q^m-nq^{m-i})=0,~\forall t,t'=0,1,\hdots,q^m-1.
	\end{split}
\end{equation}
There exist $np^{m-i}$ pairs $(x,y)$ such that $y=x+q^m-nq^{m-i}$, where $x=0,1,\hdots, nq^{m-i}-1$. Let $(x_1,x_2,\hdots,x_m)$ be the $q$-ary vector representation of $x$, i.e.,  $x=\sum_{j=1}^m x_j q^{j-1}$. Then $y=x+\tau$, implies 
\begin{equation}\label{augiff1}
	\begin{split}
	y=x+\tau=\sum_{j=1}^{m-i} x_j q^{j-1}+(x_{m-i+1}+q-n)q^{m-i}+\sum_{j=m-i+2}^{m} (x_j+q-1) q^{j-1}.
	\end{split}
\end{equation}
For $(y_1,y_2,\hdots,y_m)$ as the vector representation of $y$, 
from (\ref{augiff1}), we have $x_j=y_j$, for $j=1,2,\hdots,m-i$, $y_{m-i+1}=x_{m-i+1}+q-n$, where $x_{m-i+1}=0,1,\hdots,n-1$. Besides, $y_{j'}=x_{j'}+q-1$, where $x_{j'}=0$, for $j'=m-i+2,m-i+3,\hdots,m$. Let $S_n$ be the set of all the $nq^{m-i}$ pairs $(x,y)$ with $y-x=q^m-nq^{m-i}$, can be expressed as 
\begin{equation}\label{kuku8824}
	\begin{split}
	S_n=&\left\{(x,y): y=x+q^m-nq^{m-i},~ y_j=x_j,~y_{m-i+1}=x_{m-i+1}+q-n,~ y_{j'}=x_{j'}+q-1,\right. \\& \left. x_j=0,1,\hdots,q-1, x_{m-i+1}=0,1,\hdots,n-1,  ~\text{and}~x_{j'}=0,~\text{where}~j=1,2,\hdots,m-i,\right.\\&\left.\text{and}~j'=m-i+2,\hdots,m\right\}.
	\end{split}
\end{equation}
Then it is clear that $|S_n|=nq^{m-i}$, and 
\begin{equation}\label{key17824}
\Theta(\psi(C_t),\psi(C_{t'}))(q^m-nq^{m-i})=\sum_{(x,y)\in S_n} \xi_q^{F_{x,t}-F_{y,t'}}.
\end{equation}
Now let us assume that $n=1$, which says $\tau=q^m-q^{m-i}$, and $|S_1|=q^{m-i}$. We can partition $S_1$ into $q^{m-i-1}$ disjoint subsets with each having $q$ elements form $S$.
Let $(x^u_{\alpha},y^u_{\alpha})\in S_1$, such that $x_\alpha^u=(x_{\alpha,1},x_{\alpha,2},\hdots,x_{\alpha,m-i-1},u,0,\hdots,0)$, and 
$y_\alpha^u=(x_{\alpha,1},x_{\alpha,2},\hdots,x_{\alpha,m-i-1},u,q-1,\hdots,q-1)$, we define the following subset of $S_1$:
$$S_{1,\alpha}=\{(x_\alpha^u,y_\alpha^u):u=0,1,\hdots,q-1\},$$
where $\alpha=1,2,\hdots,q^{m-i-1}$. Note that $x_{\alpha,j}\in \mbb Z_q$, where $j=1,2,\hdots,m-i-1$, and $1\leq \alpha_1\neq \alpha_2\leq q^{m-i-1}$ implies $$(x_{{\alpha}_1,1},x_{\alpha_1,2},\hdots,x_{\alpha_1,m-i-1})\neq (x_{\alpha_2,1},x_{\alpha_2,2},\hdots,x_{\alpha_2,m-i-1}).$$    
Then $S_1=\cup_{\alpha=1}^{q^{m-i-1}} S_{1,\alpha}$, and 
\begin{equation}\label{key15824}
	\begin{split}
\Theta(\psi(C_t),\psi(C_{t'}))(q^m-q^{m-i})&=\sum_{(x,y)\in S_n} \xi_q^{F_{x,t}-F_{y,t'}}\\
&=\sum_{(x,y)\in \cup_{\alpha=1}^{q^{m-i-1}} S_{1,\alpha}} \xi_q^{F_{x,t}-F_{y,t'}}\\
&=\sum_{\alpha=1}^{q^{m-i-1}} \sum_{u=0}^{q-1} \xi_q^{F_{x_{\alpha}^u,t}-F_{y_{\alpha}^u,t'}}. 
	\end{split}
\end{equation}
Then 
\begin{equation}\nonumber
	\begin{split}
F_{x_{\alpha}^u,t}-F_{y_{\alpha}^u,t'}-(F_{x_{\alpha}^0,t}-F_{y_{\alpha}^0,t'})=(h'_{m-i}(0)-h'_{m-i}(q-1))(h_{m-i}(u)-h_{m-i}(0)),
	\end{split}
\end{equation}
which implies
\begin{equation}\nonumber
	\begin{split}
	\sum_{u=0}^{q-1} \xi_q^{F_{x_{\alpha}^u,t}-F_{y_{\alpha}^u,t'}}=\xi_q^{F_{x_{\alpha}^0,t}-F_{y_{\alpha}^0,t'}} 	\sum_{u=0}^{q-1} \xi_q^{(h'_{m-i}(0)-h'_{m-i}(q-1))(h_{m-i}(u)-h_{m-i}(0))}.
	\end{split}
\end{equation}
Then 
\begin{equation}\label{key25824}
	\begin{split}
\Theta(\psi(C_t),\psi(C_{t'}))(q^m-nq^{m-i})= \sum_{\alpha=1}^{q^{m-i-1}} \xi_q^{F_{x_{\alpha}^0,t}-F_{y_{\alpha}^0,t'}} 	\sum_{u=0}^{q-1} \xi_q^{(h'_{m-i}(0)-h'_{m-i}(q-1))(h_{m-i}(u)-h_{m-i}(0))}.
	\end{split}
\end{equation}
As in the $q$-ary vector representation of $x_{\alpha}^0$ and $y_{\alpha}^0$, for $j=1,2,\hdots,m-i-1$, $x_{\alpha,j}=y_{\alpha,j}$  and $x_{\alpha,m-i}=y_{\alpha,m-i}=0$. Also for $j'=m-i+1,\hdots,m-1$, $x_{\alpha,j'}=0$ and $y_{\alpha,j'}=q-1$,
\begin{equation}\label{key35824}
	\begin{split}
F&_{x_{\alpha}^0,t}-F_{y_{\alpha}^0,t'}\\&=h_{m-i}(0)\left(h_{m-i}'(0)-h_{m-i}'(q-1)\right)+\sum_{j=m-i+1}^{m-1}\left(h_j(0)h_j'(0)-h_j(q-1)h_j'(q-1)\right)
\\&+ \sum_{j=m-i+1}^{m-1} \left(g_j(0)-g_j(q-1)\right)+t.0-t'(q-1),
	\end{split}
\end{equation}
and it says 
\begin{equation}\nonumber
	F_{x_{\alpha_1}^0,t}-F_{y_{\alpha_1}^0,t'}=F_{x_{\alpha_2}^0,t}-F_{y_{\alpha_2}^0,t'},~\forall ~\alpha_1,~\alpha_2=1,2,\hdots, q^{m-i-1}.
\end{equation}
Assuming $F_{x_{\alpha_1}^0,t}-F_{y_{\alpha_1}^0,t'}=F_{x_{\alpha_2}^0,t}-F_{y_{\alpha_2}^0,t'}=F_{x_{\alpha}^0,t}-F_{y_{\alpha}^0,t'},~\forall ~\alpha_1,~\alpha_2=1,2,\hdots, q^{m-i-1}$ in (\ref{key25824}), we have
\begin{equation}\label{key45824}
	\begin{split}
\Theta(\psi(C_t),\psi(C_{t'}))(q^m-q^{m-i})&= \sum_{\alpha=1}^{q^{m-i-1}} \xi_q^{F_{x_{\alpha}^0,t}-F_{y_{\alpha}^0,t'}} 	\sum_{u=0}^{q-1} \xi_q^{(h'_{m-i}(0)-h'_{m-i}(q-1))(h_{m-i}(u)-h_{m-i}(0))}	\\
&=q^{m-i-1} \xi_q^{F_{x_{\alpha}^0,t}-F_{y_{\alpha}^0,t'}} 	\sum_{u=0}^{q-1} \xi_q^{(h'_{m-i}(0)-h'_{m-i}(q-1))(h_{m-i}(u)-h_{m-i}(0))}.
	\end{split}
\end{equation}
Given $\Theta(\psi(C_t),\psi(C_{t'}))(q^m-q^{m-i})=0$, (\ref{key45824}) implies that 
$h'_{m-i}(0)\neq h'_{m-i}(q-1)$, and $0<(h'_{m-i}(0)-h'_{m-i}(q-1))\mod q\leq q-1 $. 
Given the condition $h'_{m-i}(0)\neq h'_{m-i}(q-1)$, $h'_{m-i}$ can have many choices, based on which, $r=(h'_{m-i}(0)-h'_{m-i}(q-1))\mod q$ can attain all the values in $\mbb Z_q\setminus \{0\}$. 
It implies
$$\sum_{u=0}^{q-1} \xi_q^{rh_{m-i}(u)}=0~ \forall r\in \mbb Z_q\setminus \{0\},$$ which again implies $h_{m-i}$ permutes $\mbb Z_q$, where $i=1,2,\hdots,m-1$. 
For $n>1$, (\ref{key17824}) can be represented as
\begin{equation}\label{keyhat7824}
\begin{split}
\Theta(\psi(C_t),\psi(C_{t'}))(q^m-nq^{m-i})&=\sum_{(x,y)\in S_n} \xi_q^{F_{x,t}-F_{y,t'}}\\
%&=\sum_{(x,y)\in \cup_{\alpha=1}^{q^{m-i-1}} S_{1,\alpha}} \xi_q^{F_{x,t}-F_{y,t'}}\\
&=\sum_{\alpha=1}^{q^{m-i-1}}\sum_{x_{\alpha,m-i+1}=0}^{n-1} \sum_{u=0}^{q-1} \xi_q^{F_{x_{\alpha}^u,t}-F_{y_{\alpha}^u,t'}}, 
\end{split}
\end{equation}
where 
$$x_\alpha^u=(x_{\alpha,1},x_{\alpha,2},\hdots,x_{\alpha,m-i-1},u,x_{\alpha,m-i+1},0,\hdots,0),$$ and 
$$y_\alpha^u=(x_{\alpha,1},x_{\alpha,2},\hdots,x_{\alpha,m-i-1},u,x_{\alpha,m-i+1}+q-n,q-1,\hdots,q-1).$$
Then
\begin{equation}\nonumber
\begin{split}
F&_{x_{\alpha}^u,t}-F_{y_{\alpha}^u,t'}-(F_{x_{\alpha}^0,t}-F_{y_{\alpha}^0,t'})\\&=(h'_{m-i}(x_{\alpha,m-i+1})-h'_{m-i}(x_{\alpha,m-i+1}+q-n))(h_{m-i}(u)-h_{m-i}(0)),
\end{split}
\end{equation}
and
\begin{equation}\nonumber
\begin{split}
\sum_{u=0}^{q-1} \xi_q^{F_{x_{\alpha}^u,t}-F_{y_{\alpha}^u,t'}}=\xi_q^{F_{x_{\alpha}^0,t}-F_{y_{\alpha}^0,t'}} 	\sum_{u=0}^{q-1} \xi_q^{(h'_{m-i}(x_{\alpha,m-i+1})-h'_{m-i}(x_{\alpha,m-i+1}+q-n))(h_{m-i}(u)-h_{m-i}(0))}.
\end{split}
\end{equation}
Using this result, (\ref{key17824}) again can be represented as 
\begin{equation}\label{key27824}
\begin{split}
\Theta (\psi(C_t),\psi(C_{t'}))(q^m-nq^{m-i})&=\sum_{(x,y)\in S_n} \xi_q^{F_{x,t}-F_{y,t'}}\\
%&=\sum_{(x,y)\in \cup_{\alpha=1}^{q^{m-i-1}} S_{1,\alpha}} \xi_q^{F_{x,t}-F_{y,t'}}\\
&=\sum_{x_{\alpha,m-i+1}=0}^{n-1}\left(\sum_{\alpha=1}^{q^{m-i-1}} \sum_{u=0}^{q-1} \xi_q^{F_{x_{\alpha}^u,t}-F_{y_{\alpha}^u,t'}}\right), \end{split}\end{equation} \begin{equation}\nonumber
\begin{split}
&=\sum_{x_{\alpha,m-i+1}=0}^{n-1}\left(\sum_{\alpha=1}^{q^{m-i-1}} \xi_q^{F_{x_{\alpha}^0,t}-F_{y_{\alpha}^0,t'}} 	\sum_{u=0}^{q-1} \xi_q^{(h'_{m-i}(x_{\alpha,m-i+1})-h'_{m-i}(x_{\alpha,m-i+1}+q-n))(h_{m-i}(u)-h_{m-i}(0))}\right)\\
&=\sum_{j=0}^{n-1} A_j,
\end{split}
\end{equation}
where 
$A_j=\sum_{\alpha=1}^{q^{m-i-1}} \xi_q^{F_{x_{\alpha}^0,t}-F_{y_{\alpha}^0,t'}} 	\sum_{u=0}^{q-1} \xi_q^{(h'_{m-i}(j)-h'_{m-i}(j+q-n))(h_{m-i}(u)-h_{m-i}(0))}$. 
Now in the expression of $A_j$,
\begin{equation}\nonumber
\begin{split}
&F_{x_{\alpha}^0,t}-F_{y_{\alpha}^0,t'}\\&=h_{m-i}(0)\left(h'_{m-i}(j)-h'_{m-i}(j+q-n)\right)\\&+\left(h_{m-i+1}(j) h_{m-i+1}'(0)-h_{m-i+1}(j+q-n)h_{m-i+1}'(q-1)\right)\\&+\left(h_{m-i+2}(0) h_{m-i+1}'(0)-h_{m-i+1}(q-1)h_{m-i+1}'(q-1)\right)+\hdots\\&+
\left(h_{m-1}(0) h_{m-1}'(0)-h_{m-1}(q-1)h_{m-1}'(q-1)\right)+\left(g_{m-i+1}(j)-g_{m-i+1}(j+q-n)\right)\\&+\sum_{\beta=m-i+2}^m \left(g_{\beta}(0)-g_{\beta}(q-1)\right)+t\cdot 0-t'\cdot (q-1)=F_j,
\end{split}
\end{equation}
which says, for a fixed value of $j$, $F_{x_{\alpha_1}^0,t}-F_{y_{\alpha_1}^0,t'}=F_{x_{\alpha_2}^0,t}-F_{y_{\alpha_2}^0,t'}$ for all $\alpha_1,\alpha_2=1,2,\hdots,q^{m-i-1}$. Therefore, 
$$A_j=q^{m-i-1}\left(\xi_q^{F_j} 	\sum_{u=0}^{q-1} \xi_q^{(h'_{m-i}(j)-h'_{m-i}(j+q-n))(h_{m-i}(u)-h_{m-i}(0))}\right)=q^{m-i-1}B_j.$$
Then 
\begin{equation}\nonumber
	\begin{split}
	\Theta (\psi(C_t),\psi(C_{t'}))(q^m-nq^{m-i})=q^{m-i-1}\sum_{j=0}^{n-1} B_j.
	\end{split}
\end{equation} 
Now our goal is to show that for any value of $q\geq 2$, $h_{m-i}'$ also permutes $\mbb Z_q$. We assume there exist a subset $\{j_1,j_2,\hdots,j_{t}\}$ of $\{0,1,,\hdots,n-1\}$, $0\leq t\leq n$, such that $h'_{m-i}(j)-h'_{m-i}(j+q-n)=0$ for $j=j_1,j_2,\hdots,j_t$. Then 
\begin{equation}\nonumber
	\begin{split}
	\Theta (\psi(C_t),\psi(C_{t'}))(q^m-nq^{m-i})=q^{m-i}\left(\xi_q^{F_{j_1}}+\xi_q^{F_{j_2}}+\hdots+\xi_q^{F_{j_t}}\right),
	\end{split}
\end{equation} 
where $0\leq t<q-1$. Now let us assume $g_\beta=0, ~\forall \beta$, and $h_{m-i}$ is a identity permutation. Then $F_j$ can be expressed as 
$$F_j=C+j(h_{m-i+1}'(0)-h_{m-i+1}'(q-1)),$$
where $C=\left(h_{m-i+2}(0) h_{m-i+1}'(0)-h_{m-i+1}(q-1)h_{m-i+1}'(q-1)\right)+\cdots+\\
\left(h_{m-1}(0) h_{m-1}'(0)-h_{m-1}(q-1)h_{m-1}'(q-1)\right)-(q-n)h_{m-i+1}'(q-1)-t'(q-1)$, and $j=j_1,j_2,\hdots,j_t$. Therefore,
\begin{equation}\label{keyhit824}
\begin{split}
\Theta (\psi(C_t),\psi(C_{t'}))(q^m-nq^{m-i})&=q^{m-i}\left(\xi_q^{F_{j_1}}+\xi_q^{F_{j_2}}+\hdots+\xi_q^{F_{j_t}}\right)\\
&=q^{m-i}\xi_q^C\sum_{j\in\{j_1,j_2,\hdots,j_t\}} \xi_q^{j(h_{m-i+1}'(0)-h_{m-i+1}'(q-1))}.
\end{split}
\end{equation} 
Besides $0\leq t\leq n<q$, we already have $0<(h_{m-i+1}'(0)-h_{m-i+1}'(q-1))\mod q<q$, and using similar argument that we already have provided for $\tau=q^m-q^{m-i}$, it can be shown that $\Theta$ in (\ref{keyhit824}) cannot be zero unless $t=0$. Now $t=0$ implies, $h'_{m-i}(j)-h'_{m-i}(j+q-n)\neq 0,~\forall j $, which says $h'_{m-i}$ also permutes $\mbb Z_q$, where $i=1,2,\hdots,m-1$, and we complete the proof.   
\end{IEEEproof}
Below we shall extend the proposed $(q,q^m)$-CCCs to $(q^{n+1},q^m)$-CCCs by using the concept of restricting a $q$-ary functions over $n$-variables. In \cite{pater2000}, Paterson introduced the idea of restricting a GBF over some $n$ variables for extending the construction of GCPs \cite{Davis1999} to CSs of size $2^{n+1}$. Later Rathinakumar-Chaturvedi, extended Paterson's idea to construct $(2^{n+1},2^m)$-CCCs in \cite{rati}.
Let $f$ has the property that, after restricting it to $n$ variables, $x_{j_1},x_{j_2},\hdots,x_{j_n}$, the function results in a function of Hamming degree $2$ over rest of the $m-n$ variables in the following form:
\begin{equation}\label{qry_n1}
	\begin{split}
        f\arrowvert_{\mathbf{x}_J=\mathbf{c}}=\sum_{i=1}^{m-n-1} h_i(x_{\pi_c(i)})h'_i(x_{\pi_c(i+1)})+\sum_{j=1}^{m-n} g_j(x_{\pi_c(j)})+c,
	\end{split}
\end{equation} 
where $J=\{j_1,j_2,\hdots,j_n\}\subset \{1,2,\hdots,m\}$, $\mbf x_J=(x_{j_1},x_{j_2},\hdots, x_{j_n})$, $\mbf c=(c_1,c_2,\hdots,c_k)~\in\mbb Z_q^k$, $c=\sum_{i=1}^k c_i q^{k-i} $, and $\pi_c:\{1,2,\hdots,m-n\}\rightarrow \{1,2,\hdots,m\}\setminus J $ is a one-to-one mapping. The functions $h_i$, $h_i'$, and $g_j$ in (\ref{qry_n1}) represent univariate functions over $\mbb Z_q$. Corresponding to the function $ f\arrowvert_{\mathbf{x}_J=\mathbf{c}}$, we denote the complex-valued sequence
by $\psi( f\arrowvert_{\mathbf{x}_J=\mathbf{c}})$, and define as
\begin{equation}\nonumber
	\begin{split}
	\psi(f\arrowvert_{\mathbf{x}_J=\mathbf{c}})=
	\begin{cases}
	\xi_q^{f_{x}},& \textnormal{if}~\mbf x_{J}=\mbf c, \\
	0, & \textnormal{otherwise},
	\end{cases} 
	\end{split}
\end{equation}
which says, $\psi(f)=\sum_{\mbf c\in \mbb Z_q^{n}} \psi(f\arrowvert_{\mathbf{x}_J=\mathbf{c}})$. This information will be used later in the proof of our proposed construction.    
Again we assume $t_1,t_2,\hdots,t_n$ are $n$ integers from $\mbb Z_q$, and $t=\sum_{i=1}^{n+1} t_i q^{n+1-i}$. It says $0\leq t< q^{n+1}$, and $(t_1,t_2,\hdots,t_{n+1})$ is the length-$n+1$ $q$-ary vector representation of $t$. For any value of $t$, below we define $C_t$ to be a set of $q^{n+1}$ $q$-ary functions as 
\begin{equation}\nonumber
	\begin{split}
	C_t=\left\{f+\sum_{i=1}^n (d_i+t_i)x_{j_i}+d_{n+1} x_{\pi(1)}+t_{n+1} x_{\pi(m-n)}: d_i\in \mbb Z_q, i=1,2,\hdots, n+1 \right\},
	\end{split}
\end{equation}
where {$\pi=\pi_c$} if $\mbf x_J=\mbf c$.
%{$\pi(x_1)=\pi_c(x_1)=x_{\pi_c(1)}$ and $\pi(x_{m-n})=\pi_c(x_{m-n})=x_{\pi_c(m-n)}$} if {$\mbf x_J=\mbf %c$}.
Then 
\begin{equation}\nonumber
	\begin{split}
	\psi(C_t)=\left\{\psi\left(f+\sum_{i=1}^n (d_i+t_i)x_{j_i}+d_{n+1} x_{\pi(1)}+t_{n+1}x_{\pi(m-n)} \right): d_i\in \mbb Z_q, i=1,2,\hdots, n+1\right\},
	\end{split}
\end{equation}
which says $\psi(C_t)$ is a set of $q^{n+1}$ complex-valued sequences of length $q^m$ over the alphabet $\mathcal{ A}_q$. 
\ccn    
\begin{corollary}\label{cor15824}
 For any choice of the functions $g_j:\mbb Z_q\rightarrow \mbb Z_q$ , $j=1,2,\hdots,m-n$, and one-to-one mapping $\pi_c$ in (\ref{qry_n1}), the set of complex-valued codes $\mathcal{C}=\left\{\psi(C_t):t=0,1,\hdots,q^{n+1}-1 \right\}$ forms $(q^{n+1},q^m)$-CCCs iff $h_i$ and $h_i'$ permute $\mbb Z_q$, where $i=1,2,\hdots,m-n-1$. 
\end{corollary}
\begin{IEEEproof}
For any two integers $t$ and $t'$, where $0\leq t,~t'<q^{n+1}$, let us assume that $(t_1,t_2,\hdots,t_{n+1})$ and $(t_1',t_2',\hdots,t_{n+1}')$ are their $q$-ary vector representations. Then the AACF between 
$\psi(C_t)$ and $\psi(C_{t'})$ at $\tau$, where $0\leq \tau<q^m$, can be expressed as
\begin{equation}\label{labelcorr1}
	\begin{split}
	\Theta(\psi(C_t),\psi(C_{t'}))(\tau)=\sum_{\mbf d_n, d_{n+1}} \Theta\left(\psi(f+(\mbf d_n+\mbf t_n)\cdot \mbf x_J+d_{n+1} x_{\pi(1)}+t_{n+1}x_{\pi(m-n)} ),\right.\\ \left. \psi(f+(\mbf d_n+\mbf t_n')\cdot \mbf x_J+d_{n+1}x_{\pi(1)}+t_{n+1}'x_{\pi(m-n)}) \right)(\tau),
	\end{split}
\end{equation}
where $\mbf d_n=(d_1,d_2,\hdots,d_n)$  and $\mbf t_n=(t_1,t_2,\hdots,t_n)$. Again (\ref{labelcorr1}) can be represented as
\begin{equation}\label{labelcorr2}
\begin{split}
\Theta &(\psi(C_t),\psi(C_{t'}))(\tau)\\&=\sum_{\mbf c_1,~\mbf c_2\in \mbf Z_q^n}\sum_{\mbf d_n, d_{n+1}} \Theta\left(\psi\left(\left(f+(\mbf d_n+\mbf t_n)\cdot \mbf x_J+d_{n+1} x_{\pi(1)}+t_{n+1}x_{\pi(m-n)}\right)\arrowvert_{\mathbf{x}_J=\mathbf{c}_1} \right),\right.
\\&\qquad\qquad\qquad
 \left. \psi\left(\left(f+(\mbf d_n+\mbf t_n')\cdot \mbf x_J+d_{n+1} x_{\pi(1)}+t_{n+1}'x_{\pi(m-n)}\right)\arrowvert_{\mathbf{x}_J=\mathbf{c}_2}\right) \right)(\tau),\\
 &=\sum_{\mbf c_1,~\mbf c_2}\xi_q^{(\mbf t_n\cdot\mbf c_1-\mbf t_n'\mbf c_2)}\sum_{\mbf d_n}\xi_q^{\mbf d_n\cdot(\mbf c_1-\mbf c_2)}\sum_{d_{n+1}} \Theta\left(\psi\left(\left(f+d_{n+1} x_{\pi(1)}+t_{n+1}x_{\pi(m-n)}\right)\arrowvert_{\mathbf{x}_J=\mathbf{c}_1} \right),\right.
 \\&\qquad\qquad\qquad
 \left. \psi\left(\left(f+d_{n+1}x_{\pi(1)}+t_{n+1}'x_{\pi(m-n)}\right)\arrowvert_{\mathbf{x}_J=\mathbf{c}_2}\right) \right)(\tau).
\end{split}
\end{equation}
In (\ref{labelcorr2}), $$\sum_{\mbf d_n}\xi_q^{\mbf d_n\cdot(\mbf c_1-\mbf c_2)}=\begin{cases}
q^{n},& \mbf c_1=\mbf c_2,\\
0,& \mbf c_1\neq\mbf c_2.
\end{cases}$$
Then (\ref{labelcorr2}) can be represented as
\begin{equation}\label{labelcorr3}
	\begin{split}
	\Theta &(\psi(C_t),\psi(C_{t'}))(\tau)\\&=\sum_{\mbf c_1,~\mbf c_2\in \mbf Z_q^n}\xi_q^{(\mbf t_n\cdot\mbf c_1-\mbf t_n'\mbf c_2)}\sum_{\mbf d_n\in \mbb Z_q^n}\xi_q^{\mbf d_n\cdot(\mbf c_1-\mbf c_2)}\sum_{d_{n+1}\in \mbb Z_q} \Theta\left(\psi\left(\left(f+d_{n+1} x_{\pi(1)}+t_{n+1}x_{\pi(m-n)}\right)\arrowvert_{\mathbf{x}_J=\mathbf{c}_1} \right),\right.
	\\&\qquad\qquad\qquad
	\left. \psi\left(\left(f+d_{n+1}x_{\pi(1)}+t_{n+1}'x_{\pi(m-n)}\right)\arrowvert_{\mathbf{x}_J=\mathbf{c}_2}\right) \right)(\tau)\\
	&=q^{n}\sum_{\mbf c \in \mbb Z_q^n}\xi_q^{(\mbf t_n-\mbf t_n')\cdot \mbf c} \sum_{d_{n+1}\in \mbb Z_q} \Theta\left(\psi\left(\left(f+d_{n+1}x_{\pi(1)}+t_{n+1}x_{\pi(m-n)}\right)\arrowvert_{\mathbf{x}_J=\mathbf{c}} \right),\right.
	\\&\qquad\qquad\qquad
	\left. \psi\left(\left(f+d_{n+1}x_{\pi(1)}+t_{n+1}'x_{\pi(m-n)}\right)\arrowvert_{\mathbf{x}_J=\mathbf{c}}\right) \right)(\tau),
	\end{split}
\end{equation}
where $\mbf c_1=\mbf c_2=\mbf c$. For $\mbf c\in \mbb Z_q^m$, we define $\mbb Z_q^m(\mbf c)=\{\mbf x\in\mbb Z_q^m: \mbf x_J=\mbf c\}$, and $\mal N_c=\{x\in \mbf Z_{q^m}: x=\sum_{j=1}^m x_j q^{j-1},~\text{and}~\mbf x_J=\mbf c  \}$. We also define $\mathcal{ A}_\tau(\mbf c)=\{(x,y)\in  \mathcal{N}_c\times \mathcal{N}_c:y=x+\tau \}$. Now 
\begin{equation}\label{labelcorr4}
	\begin{split}
	&\sum_{d_{n+1}\in \mbb Z_q} \Theta\left(\psi\left(\left(f+d_{n+1} x_{\pi(1)}+t_{n+1} x_{\pi(m-n)} \right)\arrowvert_{\mathbf{x}_J=\mathbf{c}} \right),\right.
	\\&\qquad\qquad\qquad
	\left. \psi\left(\left(f+d_{n+1} x_{\pi(1)}+t_{n+1}'x_{\pi(m-n)}\right)\arrowvert_{\mathbf{x}_J=\mathbf{c}}\right) \right)(\tau)\\
	&=\sum_{d_{n+1}\in \mbb Z_q} \sum_{(x,y)\in \mal{A}_\tau(\mathbf{c})} \xi_q^{F_{x,t}-F_{y,t'}} \xi_q^{d_{n+1}(x_{\pi_c(1)}-y_{\pi_c(1)})},
	\end{split}
\end{equation}
 where $F_{x,t}=f_x+t_{n+1}x_{\pi_c(m-n)}$, $F_{y,t'}=f_y+t_{n+1}'y_{\pi_c(m-n)}$, $f_x=f(\mbf x)$, $f_y=f(\mbf y)$, and $\mbf y=(y_1,y_2,\hdots,y_m)$ denotes the $q$-ary vector representation of $y$. 
Now if $x_{\pi_c(1)}\neq y_{\pi_c(1)}$, it can be observed through (\ref{labelcorr4}) that $\Theta (\psi(C_t),\psi(C_{t'}))(\tau)=0$. Now we only left with the condition $x_{\pi_c(1)}=y_{\pi_c(1)}$, under which (\ref{labelcorr3}) reduces to the following: 
\begin{equation}\label{key1}
		\begin{split}
	\Theta(\psi(C_t),\psi(C_{t'}))(\tau)=q^{n+1}\sum_{\mbf c \in \mbb Z_q^n}\xi_q^{(\mbf t_n-\mbf t_n')\cdot \mbf c}\sum_{(x,y)\in \mal{A}_\tau(\mathbf{c})} \xi_q^{F_{x,t}-F_{y,t'}}.
		\end{split}
	\end{equation}
	 Now $F_{x,t}-F_{y,t'}$ can be expressed as follows:
	\begin{equation}\label{keyy}
	\begin{split}
&F_{x,t}-F_{y,t'}\\&=f_x-f_y+(t_{n+1}x_{\pi_c(m-n)}-t_{n+1}'y_{\pi_c(m-n)})\\
&=\sum_{i=1}^{m-n-1}(h_i(x_{\pi_c(i)})h'_i(x_{\pi_c(i+1)})-h_i(y_{\pi_c(i)})h'_i(y_{\pi_c(i+1)}))+\sum_{j=1}^{m-n} (g_j(x_{\pi_c(j)})-g_j(y_{\pi_c(j)}))\\& ~~~~~~~~~~+(t_{n+1}x_{\pi_c(m-n)}-t_{n+1}'{y_{\pi_c(m-n)}}).
	\end{split}
	\end{equation}
When $\tau=0$, $x=y$, which says $|\mathcal{ A}_\tau(c)|=q^{m-n}$. From (\ref{key1}) and (\ref{keyy}), 
\begin{equation}\label{key2}
\begin{split}
\Theta(\psi(C_t),\psi(C_{t'}))(0)&=q^{n+1}\sum_{\mbf c \in \mbb Z_q^k}\xi_q^{(\mbf t_n-\mbf t_n')\cdot \mbf c}\sum_{(x,x)\in \mal{A}_0(\mathbf{c})} \xi_q^{F_{x,t}-F_{x,t'}}\\
&=q^{n+1}\sum_{\mbf c \in \mbb Z_q^n}\xi_q^{(\mbf t_n-\mbf t_n')\cdot \mbf c}\sum_{(x,x)\in \mal{A}_0(\mathbf{c})} \xi_q^{(t_{n+1}-t_{n+1}')x_{\pi_c(m-n)}}\\
&=\begin{cases}
q^{m+n+1}, &\mbf t_n=\mbf t_n',~t_{n+1}=t_{n+1}',\\
0,& \text{otherwise},
\end{cases}\\
& =\begin{cases}
q^{m+n+1}, & t=t',\\
0,& t\neq t'.
\end{cases}
\end{split}
\end{equation}	
When $\tau\neq 0$, $x\neq y$, it implies $\mbf x\neq \mbf y$. Let us assume that $1<v\leq m-n$ is the smallest integer for which $x_{\pi_c(v)}\neq y_{\pi_c(v)}$. Let $(x_1,x_2,\hdots,(x_{\pi_c(v-1)}+u)\mod q,\hdots,x_m)$ and $(y_1,y_2,\hdots,(y_{\pi_c(v-1)}+u)\mod q,\hdots,y_m)$ are the $q$-ary vector representations of $x^u$ and $y^u$, respectively, where $1\leq u<q$. Besides, $y=x+\tau$ implies $y^u=x^u+\tau$, for all values of $u$. It says, if $|\mathcal{ A}_\tau(\mbf c)| \neq 0$ for some $\tau\neq 0$ and $\mbf c\in \mbb Z_k$, then $(x,y)\in \mathcal{ A}_\tau(\mbf c)$ implies $(x^u,y^u)\in \mathcal{ A}_\tau(\mbf c)$.
As $1< v\leq m-n$, $x^u$ differs with $x$ only at the $\pi(v-1)$th component,
\begin{equation}\label{key3}
	\begin{split}
	F_{x^u,t}-F_{x,t} &=f_{x^u}-f_x\\&= h_{v-2}(x_{\pi_c(v-2)})\left(h'_{v-2}((x_{\pi_c(v-1)}+u)\mod q)-h'_{v-2}(x_{\pi_c(v-1)})   \right)+\\&
	~~~h'_{v-1}(x_{\pi_c(v)})\left( h_{v-1}((x_{\pi_c(v-1)}+u)\mod q)-h_{v-1}(x_{\pi_c(v-1)})  \right)+\\&~~~\left(g_{v-1}((x_{\pi_c(v-1)}+u)\mod q)-g_{v-1}(x_{\pi_c(v-1)}) \right).
	\end{split}
\end{equation} 
As $\pi_c(x_i)=\pi_c(y_i)$ for $i=1,2,\hdots,v-1$, similarly as (\ref{key3}), we obtain $F_{y^u,t'}-F_{y,t'}$ as follows:
\begin{equation}\label{key4}
	\begin{split}
	F_{y^u,t'}-F_{y,t'} &=f_{y^u}-f_y\\&= h_{v-2}(x_{\pi_c(v-2)})\left(h'_{v-2}((x_{\pi_c(v-1)}+u)\mod q)-h'_{v-2}(x_{\pi_c(v-1)})   \right)+\\&
	~~~h'_{v-1}(y_{\pi_c(v)})\left( h_{v-1}((x_{\pi_c(v-1)}+u)\mod q)-h_{v-1}(x_{\pi_c(v-1)})  \right)+\\&~~~\left(g_{v-1}((x_{\pi_c(v-1)}+u)\mod q)-g_{v-1}(x_{\pi_c(v-1)}) \right).
	\end{split}
\end{equation}
Then 
\begin{equation}\nonumber
	\begin{split}
	&F_{x^u,t}-F_{y^u,t'}-(F_{x,t}-F_{y,t'})\\&=(F_{x^u,t}-F_{x,t})-(F_{y^u,t'}-F_{y,t'})\\&
	=\left(h'_{v-1}(x_{\pi_c(v)})-h'_{v-1}(y_{\pi_c(v)})\right)\left( h_{v-1}((x_{\pi_c(v-1)}+u)\mod q)-h_{v-1}(x_{\pi_c(v-1)})  \right)
	\end{split}
\end{equation}
Hence,
\begin{equation}\label{key5}
	\begin{split}
\sum_{u=1}^{q-1}\xi_q^{F_{x^u,t}-F_{y^u,t'}-(F_{x,t}-F_{y,t'})}=\sum_{u=1}^{q-1}\xi_q^{rp(u)}=\left(\sum_{u=0}^{q-1}\xi_q^{rp(u)}\right)-1,
	\end{split}
\end{equation}
equivalently, we can represent this as 
\begin{equation}\nonumber
\begin{split}
\xi_q^{F_{x,t}-F_{y,t'}}+\sum_{u=1}^{q-1}\xi_q^{F_{x^u,t}-F_{y^u,t'}}=\sum_{u=1}^{q-1}\xi_q^{rp(u)}=\left(\sum_{u=0}^{q-1}\xi_q^{rp(u)}\right),
\end{split}
\end{equation}
where $r=\left(h'_{v-1}(x_{\pi_c(v)})-h'_{v-1}(y_{\pi_c(v)})\right)$ and $p(u)=\left( h_{v-1}((x_{\pi_c(v-1)}+u)\mod q)-h_{v-1}(x_{\pi_c(v-1)})  \right)$. 
\setcounter{case}{0}
\noindent \textit{\textbf{\textbf{Sufficiency:}}}
Let us assume $h_i$ and $h_i'$ permute $\mbb Z_q$, where $i=1,2,\hdots,m-n-1$. It implies $r\neq 0$ and $p(u)$ permutes $\mbb Z_q$, and it again implies
$$\xi_q^{F_{x,t}-F_{y,t'}}+\sum_{u=1}^{q-1}\xi_q^{F_{x^u,t}-F_{y^u,t'}}=0.$$ 
Then by combining all the above results, we reach the following result: 
\begin{equation}\label{8824newcond2_12724}
\Theta(\psi(C_t),\psi(C_{t'}))(\tau)=0~\text{if}~h_{v-1},~ h_{v-1}'~\text{permute}~\mbb Z_q,
\end{equation}
where $v=2,3,\hdots,m-n$. Now for $\tau<0$, and applying similar logic as above, we can also reach to the same result as in (\ref{8824newcond2_12724}).

\noindent \textit{\textbf{Necessity:}}
Let us assume that, for any choice of the functions $g_j:\mbb Z_q\rightarrow \mbb Z_q$, $j=1,2,\hdots,m-n$, and one-to-one mappings $\pi_c$ in (\ref{qry_n1}), the set of complex-valued codes $\mathcal{C}=\{\psi(C_t):t=0,1,\hdots,q^{n+1}-1 \}$ forms $(q^{n+1},q^{m})$-CCCs. Using this property, our goal is to determine the type of the functions $h_i$ and $h_i'$, where $i=1,2,\hdots,m-n-1$. Let us assume $\tau=q^{k}-rq^{k-i}$, where $k=m-n$, $1\leq r\leq q-1$, and $i=1,2,\hdots,k-1$. As $\mathcal{C}$ forms CCCs, 
\begin{equation}\label{8824key1}
\begin{split}
\Theta(\psi(C_t),\psi(C_{t'}))(\tau)=q^{n+1}\sum_{\mbf c \in \mbb Z_q^n}\xi_q^{(\mbf t_n-\mbf t_n')\cdot \mbf c}\sum_{(x,y)\in \mal{A}_\tau(\mathbf{c})} \xi_q^{F_{x,t}-F_{y,t'}}=0.
\end{split}
\end{equation}
To reach our goal, i.e., to determine the properties of $h_i$ and $h_i'$, it is enough to take care of the term $\sum_{(x,y)\in \mal{A}_\tau(\mathbf{c})} \xi_q^{F_{x,t}-F_{y,t'}}$ in (\ref{8824key1}). To do so, let us assume $\mbf x_J=(x_{k+1},x_{k+2},\hdots,x_{m})$ and $\pi_c(i')=i'$~$\forall c$, where $i'=1,2,\hdots,k$.
Now let us focus on $\mal{A}_\tau(\mbf c)$. Let $(x,y)\in \mal{A}_\tau(\mbf c)$, then it is to be noted that in the $q$-ary vector representations of $x$ and $y$, $\mbf x_J=\mbf c$ and $\mbf y_J=\mbf c$. Now $y=x+\tau$ implies $$y=\sum_{j=1}^{k-1}x_jq^{j-1}+(x_{k-i+1}+q-r)q^{k-i}+(x_{k-i+2}+q-1)q^{k-i+1}+\cdots+(x_k+q-1)q^{k-1}+c',$$
where $c'=c_1q^k+c_2q^{k+1}+\cdots+c_nq^{m-1}$. Then 
\begin{equation}\label{kuku18824}
\begin{split}
\mal{A}_\tau(\mbf c)=&\left\{(x,y): y=x+q^k-rq^{k-i},~ y_j=x_j,~y_{k-i+1}=x_{k-i+1}+q-r,~ y_{j'}=x_{j'}+q-1,\right. \\& \left. x_j=0,1,\hdots,q-1, x_{k-i+1}=0,1,\hdots,r-1,  ~\text{and}~x_{j'}=0,~\text{where}~j=1,2,\hdots,k-i,\right.\\&\left.~j'=k-i+2,\hdots,k-1,~\text{and}~\mbf x_J=\mbf c=\mbf y_J\right\}.
\end{split}
\end{equation} 
Then comparing $S_n$ in (\ref{kuku8824}), and $\mal{A}_\tau(\mbf c)$ in (\ref{kuku18824}), it can be observed that they basically follow similar properties. Therefore we omit the proof as the rest of the proof will follow the similar logic as presented in the \textbf{\textit{Necessity}} part in the proof of \textbf{Theorem \ref{th824}}. So we conclude this proof here that, similarly to the proof of \textbf{Theorem \ref{th824}},
it can be shown that the code set $\mal{C}$ forms CCCs implies $h_i$ and $h_i'$ permute $\mbb Z_q$, where $i=1,2,\hdots,m-n-1$. 
\end{IEEEproof}
Now we consider the expression \( L = p_1^{m_1} p_2^{m_2} \), where \( p_1 \) and \( p_2 \) are distinct primes, and \( m_1 \neq m_2 \). To generate CCCs of length \( L \) using \textbf{Corollary \ref{cor15824}}, there is only one possible selection for the parameters \( q \), \( n \), and \( m \), which is \( q = L \), \( n = 0 \), and \( m = 1 \). This leads to the generation of \((L, L)\)-CCCs.
Consequently, the alphabet size in this scenario is also \( L \), i.e., when $m_1\neq m_2$, using \textbf{Corollary \ref{cor15824}}, we can only obtain those CCCs for which the alphabet size is as large as sequence length that may not be feasible for its practical application.

Moreover, it has been demonstrated in \cite[Th. 1]{jinla} that given a \((K_1, L_1)\)-CCCs \( \mathcal{C} = \{C_u : u = 1, 2, \dots, K_1\} \) and a \((K_2, L_2)\)-CCCs \( \mathcal{D} = \{D_v : v = 1, 2, \dots, K_2\} \), the set \( \{C_u \otimes D_v : u = 1, 2, \dots, K_1, v = 1, 2, \dots, K_2 \} \) forms a \((K_1 K_2, L_1 L_2)\)-CCCs. Applying \textbf{Corollary \ref{cor15824}} again, we can obtain \((p_1^{n_1+1}, p_1^{m_1})\)-CCCs over the alphabet \( \mal A_{p_1} \) and \((p_2^{n_2+1}, p_2^{m_2})\)-CCCs over the alphabet \( \mal A_{p_2} \). By setting \( K_1 = p_1^{m_1} \), \( L_1 = p_1^{m_1} \), \( K_2 = p_2^{n_2+1} \), and \( L_2 = p_2^{m_2} \), and using the aforementioned result of \cite{jinla}, we can construct \((p_1^{n_1+1} p_2^{n_2+1}, L)\)-CCCs over an alphabet of size \( p_1 p_2 < q \).

This demonstrates that in many cases, particularly when \( m_1 \neq m_2 \), \textbf{Corollary \ref{cor15824}} can not maintain a small alphabet size, and also the set size \( p_1^{n_1+1} p_2^{n_2+1} \leq L \) as \( 0 \leq n_1 \leq m_1 - 1 \) and \( 0 \leq n_2 \leq m_2 - 1 \). According to existing literature, having a CS with a smaller set size is advantageous for many practical applications, such as PMEPR reduction in OFDM systems. This limitation of \textbf{Corollary \ref{cor15824}} motivates us for the further investigation on the functional forms of \((K, L)\)-CCCs, where \( K = \prod_{i=1}^k p_i^{n_i+1} \), \( L = \prod_{i=1}^k p_i^{m_i} \), and the alphabet size can be maintained as low as \( \prod_{i=1}^k p_i \), with \( p_1, p_2, \dots, p_k \) being distinct primes.
\section{Construction CCCs over Small Alphabet Using $q$-ary function}\label{sectn4}
In the previous section, we have introduced construction of $(q^{n+1},q^m)$-CCCs over the alphabet $\mathcal{ A}_q$ using $q$-ary functions with domain $\mbb Z_q^m$. In this section, we introduce the use of another type of $q$-ary functions with the domain $\mbf Z_{p_1}^{m_1}\times \mbf Z_{p_2}^{m_2}\times \cdots\times \mbf Z_{p_k}^{m_k}$, where  $q=\prod_{i=1}^k p_i$.  
\subsection{Functions and Traditional Sequences}\label{seca16824}
In this section, we define necessary algebraic operations on the set $\mal{V}_L$ to introduce $q$-ary functions while maintaining $\mal{V}_L$ as the domain. 

Let $\mathbf{x}  = (\bx_1, , \dots, \bx_k) = (x_1, \dots, x_{m_1}, \dots, x_{m_1+\cdots+m_{k-1}+1}, \dots, x_{m})\in \mal V_L$, where 
$\mbf x_1=(x_1,x_2,\hdots,x_{m_1})~\in \mbf Z_{p_1}^{m_1}$, \textnormal{and}
$\bx_i = (x_{m_1+\cdots+m_{i-1}+1}, \dots, x_{m_1+\hdots+m_i}) \in \bZ_{p_i}^{m_i},~2\leq i\leq k$. 
Let us also assume that $\bx$ is the vectorial representation of an integer $x \in \bZ_L$ and $x$ is the integer representation of a vector $\bx \in \mal{V}_L$, which are connected via the following equality
\begin{equation}\label{18924connection1}
x=\sigma_1\Delta_1+\sigma_2\Delta_2+\hdots+\sigma_k\Delta_k,
\end{equation}
%$$x=\rho^{-1}(\mbf x)=\sigma_1\Delta_1+\sigma_2\Delta_2+\hdots+\sigma_k\Delta_k,$$
where  $\Delta_i = \frac{L}{L_iL_{i+1}\cdots L_{k}}$, $L_i=p_i^{m_i}$, $i=1,2,\hdots,k$. 
Also, 
\begin{equation}
\begin{split}
\sigma_1&=\sum_{j=1}^{m_1}x_jp_1^{j-1},~\textnormal{and}\\
\sigma_i&=\sum_{j=1}^{m_i} x_{m_1+\cdots+m_{i-1}+j} p_i^{j-1},~i=2,3,\hdots,k.
\end{split}
\end{equation}
The correspondence in $\mbf x$ and $x$, as outlined in (\ref{18924connection1}), is explained in greater detail in Appendix \ref{appendix1}. \ccn
Now we consider functions from $\calV_L$ to $\bbZ_q$ in $m$ variables $x_{1},x_2,\hdots,x_m$,
where the arithmetic operations among these variables and coefficients are taken modulo $q$. 
For an element $\be = (\be_1,\dots, \be_k)=(e_1,e_2\hdots,e_m) \in \calV_L$, we define a monomial $\mbf{x}^\mbf{e}$ over the $m$ variables as
\[
\bx^\be = \bx_1^{\be_1} \cdots \bx_k^{\be_k} = \prod_{j=1}^{m} x_{j}^{e_{j}},
\] where $\mbf e_1=(e_1,e_2,\hdots,e_{m_1})~\in \mbf Z_{p_1}^{m_1}$, and $\be_i = (e_{m_1+\cdots+m_{i-1}+1}, \dots, e_{m_1+\hdots+m_i}) \in \bZ_{p_i}^{m_i},~2\leq i\leq k$. Also, by convention, we assume $0^0=0$ and $x^0 =1$ for $x\neq 0$.  
A linear combination of the monomials $\mbf{x}^\mbf{e}$, with $\mbf e \in \mal V_L$ and $\mathbb{Z}_q$-valued coefficients gives
a $q$-ary function $f:\, \calV_L \rightarrow \bbZ_q$ which can be expressed as
\[
f(\bx) = \sum_{\be \in \calV_L} c_{\be} \bx^\be = \sum_{\be \in \calV_L} c_{\be} \bx_1^{\be_1} \cdots \bx_k^{\be_k}.
\] 
For $k=1$ and $p_1=2$, the function $f$ reduces to a Boolean function \cite{Davis1999}, and for $p_1>2$, it reduces to a extended Boolean function \cite{shen2023}. 
We define \textit{Hamming degree} of the above function $f$ as 
\[
\deg_H(f) = \max\{\wt(\be)\,:\, c_{\be} \neq 0 \},
\] where $\wt(\be)$ is the Hamming weight of $\be = (e_{1},e_2,\hdots,e_m)$. 
%From the definition, it is clear that any $q$-ary function $f$ has $\deg_H(f) \leq \deg(f)$.
%We will be concerned with a special class of $q$-ary functions from $\calV_L$ to $\bbZ_q$ which have Hamming degree at most $2$, namely,
A $q$-ary functions from $\calV_L$ to $\bbZ_q$ of Hamming degree at most $r$ can be \ccr uniquely \ccn expressed as
\[
f(\bx) = \sum_{\be \in \calV_L, \wt(\be)\leq r} c_{\be} \bx^\be.
\] 
Denote by $\mho_{L,r}$ the set of all the above $q$-ary functions with Hamming degree at most $r$. Then $\mho_{L,r}$ can be expressed as
\begin{equation}\nonumber
	\mho_{L,r}=\left\{\sum_{\be \in \calV_L, \wt(\be)\leq r} c_{\be} \bx^\be: c_\be\in \mathbb{Z}_q \right\}.
\end{equation} 
%%%%%%%%%%%%%%%%%%%%%%%%%%%%%%%%
%%%%%%%%%%Example 1%%%%%%%%%%%%%

\begin{example}\label{ankita008}
	Let us assume $k=2$, $p_1=2$, $p_2=3$, $m_1=3$, $m_2=2$, $q=p_1 p_2=6$, and so $L=p_1^{m_1}p_2^{m_2}=72$.   
	%Then $\mbf{x}=(x_{1,1},x_{1,2},x_{1,3},x_{2,1},x_{2,2})\in \calV_{72}$.  t us define $f:\calV_{72}\rightarrow \mathbb{Z}_6$, where $\calV_{72}$ is presented in 
	Table \ref{dtable} containing the $q$-ary vector representation $\mbf{x}=(x_{1},x_{2},x_{3},x_{4},x_{5})\in \calV_{72}$ of $x$, where $x\in \mbf{Z}_{72}$, and Table \ref{ankita700} contains all the monomials of Hamming weight at most $2$ over the $5$ variables $x_{1},x_{2},x_{3},x_{4},x_{5}$. There are $27$ monomials and therefore $\mho_{72,2}$ contains $6^{27}$ functions which can be obtained by taking a linear combination of the monomials in Table\ref{ankita700} with coefficients from $\mathbb{Z}_6$.     
\begin{table}[H]
	\centering
	\caption{The Set $\calV_{72}$ for $p_1=2$, $p_2=3$, $m_1=3$, and $m_2=2$.}\label{dtable}
	\begin{tabular}{|l|l|l|l|l|l|l|l|l|l|l|l|}
		\hline
		$x$&$\mbf{x}$&$x$&$\mbf{x}$&$x$&$\mbf{x}$&$x$&$\mbf{x}$&$x$&$\mbf{x}$&$x$&$\mbf{x}$\\ \hline                                       
		\begin{tabular}[c]{@{}l@{}}\ggr $0$\\ \ggr $1$\\ \ggr $2$\\ \ggr $3$\\\ggr $4$\\ \ggr \ggr $5$\\\ggr $6$\\ $7$\\\ggr $8$\\\ggr $9$\\\ggr $10$\\ $11$\ccn \end{tabular} & \begin{tabular}[c]{@{}l@{}}$0     0     0     0     0$ \\  $1     0     0     0     0$\\  $0     1     0     0     0$\\  $1     1     0     0     0$\\  $0     0     1     0     0$\\  $ 1     0     1     0     0$\\  $0     1     1     0     0$\\  $1     1     1     0     0$\\  $0     0     0     1     0$\\  $1     0     0     1     0$\\  $0     1     0     1     0$\\  $1     1     0     1     0$\end{tabular} & \begin{tabular}[c]{@{}l@{}}\ggr $12$\\ $13$\\ $14$\\ $15$\\\ggr $16$\\\ggr $17$\\\ggr $18$\\ $19$\\\ggr $20$\\ $21$\\ $22$\\  $23$\end{tabular} & \begin{tabular}[c]{@{}l@{}}$   0     0     1     1     0$\\     $ 1     0     1     1     0$\\     $ 0     1     1     1     0$\\      $1     1     1     1     0$\\      $0     0     0     2     0$\\      $1     0     0     2     0$\\      $0     1     0     2     0$\\      $1     1     0     2     0$\\      $0     0     1     2     0$\\      $1     0     1     2     0$\\      $0     1     1     2     0$\\      $1     1     1     2     0$\end{tabular} & \begin{tabular}[c]{@{}l@{}}\ggr $24$\\\ggr $25$\\\ggr $26$\\ $27$\\\ggr $28$\\  $29$\\ $30$\\ $31$\\\ggr $32$\\ $33$\\ $34$\\ $35$\end{tabular} & \begin{tabular}[c]{@{}l@{}}$ 0     0     0     0     1$\\     $ 1     0     0     0     1$\\     $ 0     1     0     0     1$\\     $ 1     1     0     0     1$\\     $ 0     0     1     0     1$\\     $ 1     0     1     0     1$\\     $ 0     1     1     0     1$\\     $ 1     1     1     0     1$\\     $ 0     0     0     1     1$\\     $ 1     0     0     1     1$\\     $ 0     1     0     1     1$\\     $ 1     1     0     1     1$\end{tabular} & \begin{tabular}[c]{@{}l@{}} $36$\\  $37$\\ $38$\\ $39$\\\ggr $40$\\ $41$\\ $42$\\ $43$\\ $44$\\  $45$\\ $46$\\ $47$\end{tabular} & \begin{tabular}[c]{@{}l@{}}$ 0     0     1     1     1$\\     $1     0     1     1     1$\\     $ 0     1     1     1     1$\\     $ 1     1     1     1     1$\\     $ 0     0     0     2     1$\\     $ 1     0     0     2     1$\\     $ 0     1     0     2     1$\\     $ 1     1     0     2     1$\\     $ 0     0     1     2     1$\\     $ 1     0     1     2     1$\\     $ 0     1     1     2     1$\\     $ 1     1     1     2     1$\end{tabular} & \begin{tabular}[c]{@{}l@{}}\ggr $48$\\\ggr $49$\\\ggr $50$\\ $51$\\\ggr $52$\\  $53$\\ $54$\\ $55$\\\ggr $56$\\ $57$\\ $58$\\ $59$\end{tabular} & \begin{tabular}[c]{@{}l@{}}
			$0     0     0     0     2$\\     $ 1     0     0     0     2$\\      $0     1     0     0     2$\\      $1     1     0     0     2$\\      $0     0     1     0     2$\\      $1     0     1     0     2$\\      $0     1     1     0     2$\\      $1     1     1     0     2$\\      $0     0     0     1     2$\\      $1     0     0     1     2$\\      $0     1     0     1     2$\\      $1     1     0     1     2$\end{tabular} & \begin{tabular}[c]{@{}l@{}}$60$\\ $61$\\ $62$\\ $63$\\\ggr $64$\\ $65$\\ $66$\\ $67$\\ $68$\\ $69$\\ $70$\\ $71$\end{tabular} & \begin{tabular}[c]{@{}l@{}}$0     0     1     1     2$\\      $1     0     1     1     2$\\      $0     1     1     1     2$\\      $1     1     1     1     2$\\      $0     0    0     2     2$\\      $1     0     0     2     2$\\      $0     1     0     2     2$\\      $1     1     0     2     2$\\      $0     0     1     2     2$\\      $1     0     1     2     2$\\      $0     1     1     2     2$\\      $1     1     1     2     2$\end{tabular} \\ \hline
	\end{tabular}
\end{table}

	\begin{table}[H]
		\centering
		\caption{Set of Monomials of Hamming weight at most $2$.}\label{ankita700}
		\begin{tabular}{|l|l|l|l|l|l|l|l|l|l|l|l|l|}
			\hline
			$\mbf{e}\in \calV_{72}$&$\mbf{x}^\mbf{e}$&   $\mbf{e}\in \calV_{72}$  &   $\mbf{x}^\mbf{e}$  &  $\mbf{e}\in \calV_{72}$  & $\mbf{x}^\mbf{e}$   \\ \hline                                       
			\begin{tabular}[c]{@{}l@{}}$0     0     0     0     0$ \\  $0     0     0     0     1$\\  $0     0     0     0     2$\\  $0     0     0     1     0$\\  $0     0     0     1     1$\\  $ 0     0     0     1     2$\\  $0     0     0     2     0$\\  $0     0     0     2     1$\\  $0     0     0     2     2$\end{tabular} & 
			\begin{tabular}[c]{@{}l@{}} $1$\\ $x_{5}$\\ $x_{5}^2$\\ $x_{4}$\\ $x_{4}x_{5}$\\ $x_{4}x_{5}^2$\\ $x_{4}^2$\\ $x_{4}^2x_{5}$\\ $x_{4}^2 x_{5}^2$\end{tabular} 
			& 
			\begin{tabular}[c]{@{}l@{}} $0     0     1     0     0$\\  $0     0     1     0     1$\\  $0     0     1     0     2$\\$   0     0     1     1     0$\\  $0     0     1     2     0$\\           $0     1     0     0     0$\\      $0     1     0     0     1$\\      $0     1     0     0     2$\\      $0     1     0     1     0$ \end{tabular} & \begin{tabular}[c]{@{}l@{}} $x_{3}$\\ $x_{3}x_{5}$\\  $x_{3}x_{5}^2$\\ $x_{3}x_{4}$\\ $x_{3}x_{4}^2$\\ $x_{2}$\\ $x_{2}x_{5}$\\ $x_{2}x_{5}^2$\\ $x_{2}x_{4}$\end{tabular} &
			\begin{tabular}[c]{@{}l@{}}$ 0     1     0     2     0$\\ $ 0     1     1     0     0$\\$ 1     0     0     0     0$\\     $1     0     0     0     1$\\     $ 1     0     0     0     2$\\     $ 1     0     0     1     0$\\ $ 1     0     0     2     0$\\ $ 1     0     1     0     0$\\ $1     1     0     0     0$ \end{tabular} 
			&
			\begin{tabular}[c]{@{}l@{}}
				$x_{2}x_{4}^2$\\ $x_{2}x_{3}$ \\ $x_{1}$\\ $x_{1} x_{5}$\\ $x_{1}x_{5}^2$\\ $x_{1}x_{4}$\\ $x_{1}x_{4}^2$\\ $x_{1}x_{3}$\\$x_{1}x_{2}$
			\end{tabular} \\
			\hline
		\end{tabular}
	\end{table}
\end{example}
%%%%%%%%%%%%%%%%%%%%%%%%%%%%%%%%%%%%%%%%%%%%%%%
%%%%%%%%%%%Example End Here%%%%%%%%%%%%%%%%%%%%
For a $q$-ary function $f:\calV_L\rightarrow \mathbb{Z}_q$, we define the $\mbb{Z}_q$-valued sequence of length $L$ as follows:
\begin{equation}\label{ajkit1}
	\begin{split}
	\eta(f)=(f_0,f_1,\hdots,f_{L-1}),
	\end{split}
\end{equation}
where $f_x=f(\mbf{x})$ and $x\in\mbf{Z}_L$. 
We define the complex-valued sequence over the alphabet $\mbb{Z}_q$ as 
\begin{equation}\label{seq007}
	\psi(f)=(\om^{f_0},\om^{f_1},\hdots,\om^{f_{L-1}}).
\end{equation}
As all the components of $\psi(f)$ are non-zero, $\psi(f)$ is a traditional sequence over the alphabet $\mathbb{Z}_q$, and (\ref{seq007}) establish a relationship between $f$ with a length-$L$ traditional complex-valued sequence $\psi(f)$. 
\begin{example}\label{exmo1} 
In Example \ref{ankita008}, let us consider the following $6$-ary function $f\in \mho_{72,2}$: 
\begin{equation}\label{ankby1}
f(x_{1},x_{2},x_{3},x_{4},x_{5})=
2x_{1}x_{2}+4x_{2}x_{3}+x_{2}x_{4}+x_{2}x_{5}+3x_{1}x_{3}+2x_{4}x_{5}
+x_{2}+2.
\end{equation}
From (\ref{ajkit1}), the $\mathbb{Z}_6$-valued sequence corresponding to $f$ is of length, $L=2^3 3^2=72$, and appears as 
$$\eta(f)=(2     2     3     5     2     5     1     0     2     2     4     0     2     5     2     1     2     2     5     1     2     5     3     2     2     2     4     0     2     5     2     1     4     4     1     3     4     1     5     4     0     0     4     0
0     3     2     1     2     2     5     1     2     5     3     2     0     0     4     0     0     3     2     1     4     4     3     5     4     1     1     0),$$
and the complex-valued sequence corresponding to $f$ is 
\begin{equation}\nonumber
\begin{split}
\psi(f)=(     \zeta_6^2     \zeta_6^2     \zeta_6^3     \zeta_6^5     \zeta_6^2     \zeta_6^5     \zeta_6^1     \zeta_6^0     \zeta_6^2     \zeta_6^2     \zeta_6^4     \zeta_6^0     \zeta_6^2     \zeta_6^5     \zeta_6^2     \zeta_6^1     \zeta_6^2     \zeta_6^2     \zeta_6^5     \zeta_6^1     \zeta_6^2     \zeta_6^5     \zeta_6^3     \zeta_6^2     \zeta_6^2     \zeta_6^2     \zeta_6^4     \zeta_6^0     \zeta_6^2     \zeta_6^5     \zeta_6^2     \zeta_6^1     \zeta_6^4     \zeta_6^4     \zeta_6^1     \zeta_6^3     \zeta_6^4     \zeta_6^1\\     \zeta_6^5     \zeta_6^4     \zeta_6^0     \zeta_6^0     \zeta_6^4     \zeta_6^0
\zeta_6^0     \zeta_6^3     \zeta_6^2     \zeta_6^1     \zeta_6^2    \zeta_6^2     \zeta_6^5     \zeta_6^1     \zeta_6^2     \zeta_6^5     \zeta_6^3     \zeta_6^2     \zeta_6^0     \zeta_6^0     \zeta_6^4     \zeta_6^0     \zeta_6^0     \zeta_6^3     \zeta_6^2     \zeta_6^1     \zeta_6^4    \zeta_6^ 4     \zeta_6^3     \zeta_6^5     \zeta_6^4     \zeta_6^1     \zeta_6^1     \zeta_6^0).
\end{split}
\end{equation}	
\end{example}
In the below section, we define details about sequences corresponding to the restriction of $q$-ary functions.\ccn
\subsection{Restriction of Functions and Their Complex-Valued Sequences}\label{sectn5}
Similarly as $\calV_L$, we define another set $\calV_{L'}$ containing $L'$ vectors of length $\sum_{i=1}^k n_i$, representing the integers $0,1,\hdots,\prod_{i=1}^kp_i^{n_i}-1$, where $L'=\prod_{i=1}^k L'_{i}$, and  $L'_i=p_i^{n_i}$. Then $\calV_{L'}$
can be expressed as
$$\calV_{L'}=\left\{\mathbf{c}=(\mathbf{c}_1,\mathbf{c}_2,\hdots,\mathbf{c}_k):  \mathbf{c}_i\in\mathbf{Z}_{p_i}^{n_i},~i=1,2,\hdots,k\right\}=\mathbf{Z}_{p_1}^{n_1}\times \mathbf{Z}_{p_2}^{n_2}\times \cdots\times\mathbf{Z}_{p_k}^{n_k},$$
where $\mbf c_i=(c_{i,1},c_{i,2},\hdots, c_{i,n_i})$, $i=1,2,\hdots,k$.
%where $n=\sum_{i=1}^k n_i$.
When $n_i=m_i$, $\calV_{L'}=\calV_{L}$.
Let $J=\{J_1,J_2,\hdots,J_k\}$ $\subset \{1,2,\hdots,m\}$, $J_i=\{j_{i,1},j_{i,2},\hdots,j_{i,n_i}\}\subset \{\sum_{\alpha=1}^{i-1}m_i +1,\hdots,\sum_{\alpha=1}^{i} m_{\alpha}\}$, and $\mathbf{x}_{J}=({\mathbf{x}_1}_{J_1},{\mathbf{x}_2}_{J_2},$ $\hdots,{\mathbf{x}_k}_{J_k})$, where ${\mathbf{x}_i}_{J_i}=({x}_{j_{i,1}},{x}_{j_{i,2}},\hdots,{x}_{j_{i,n_i})}$, $i=1,2,\hdots,k$. 

When $\mathbf{x}_{J}=\mbf{c}$ for some $\mbf{c}\in \calV_{L'}$ in $f$, we call it a restriction of $f$ and we denote it by
$f\arrowvert_{\mbf{x}_J=\mbf{c}}$. To define sequences corresponding to a $q$-ary function having restriction on some variables, let us first define the following sets
$$\calV_{L}^c=\{\mbf{x}\in \calV_{L}: {\mathbf{x}_1}_{J_1}=\mbf{c}_1,{\mathbf{x}_2}_{J_2}=\mbf{c}_2,\hdots,{\mathbf{x}_k}_{J_k}=\mbf{c}_k\}\subset \calV_{L},$$ and $$\mal{N}_c=\{x\in \mbf{Z}_L:x=\rho^{-1}(\mbf{x}), \mbf{x}\in\calV_{L}^c\}\subset \mbf{Z}_L,$$
where $c=c_1\Delta_1'+c_2\Delta_2'+\hdots+c_k\Delta_k'\in \mbf{Z}_{L'}$, $c_i=\sum_{j=1}^{n_i} c_{i,j} p_i^{j-1}\in \mbf Z_{L_i'}$, and $\Delta_i'=\frac{L'}{L_i'L_{i+1}'\cdots L_k'}$. 

In other word $c$ is the integer representation of $\mbf c$.
Then
\begin{equation}\label{nul1}
	\begin{split}
\calV_{L}=\displaystyle{\cup_{c\in \mbf Z_{L'}}} {\calV_{L}^c},~
\text{and}~
\mathbf{Z}_L=\displaystyle{\cup_{c\in \mbf Z_{L'}}} {\mal{N}_c}.
	\end{split}
\end{equation}
Also for any $c\neq c'$ in $\mbf Z_{L'}$,
\begin{equation}\label{nul2}
	\begin{split}
	\calV_{L}^c\cap \calV_{L}^{c'}=\emptyset, ~\text{and}~\mal{N}_{c}\cap \mal{N}_{c'}=\emptyset.
	\end{split}
\end{equation}
We define a length-$L$ complex-valued sequence corresponding to $f\arrowvert_{\mbf{x}_J=\mbf{c}}$ as follows:
\begin{equation}\label{lbl1}
\psi(f\arrowvert_{\mbf{x}_J=\mbf{c}})=
\begin{cases}
\zeta_q^{f_x}, & x\in \mal{N}_{c},\\
0,& \text{otherwise},
\end{cases}
\end{equation}
where $f_x=f(\mbf{x})$, and $\mbf{x}\in \calV_{L}^c$. 
From (\ref{nul1}), (\ref{nul2}), and (\ref{lbl1}), we have
\begin{equation}\label{slit1}
	\begin{split}
\psi(f)=\sum_{\mbf{c}\in\mal{V}_{L'}} \psi(f\arrowvert_{\mbf{x}_J=\mbf{c}}),
	\end{split}
\end{equation}
which says the complex-valued sequence corresponding to a $q$-ary function $f$ can be expressed as the sum of $\prod_{i=1}^k p_i^{n_i}$
length-$L$ complex-valued sequences.
In the later presentation, whenever we use the notation $\tau$, we consider it as non-negative integer. 
For two $\mathbb{Z}_q$-valued functions $f$ and $g$, from (\ref{lbl1}), the ACCF between $\psi(f\arrowvert_{\mbf{x}_J=\mbf{c}})$ and $\psi(g\arrowvert_{\mbf{x}_J=\mbf{c}})$ at $\tau$ can be expressed as 
\begin{equation}\label{lbl3}
\Theta(\psi(f\arrowvert_{\mbf{x}_J=\mbf{c}}),\psi(g\arrowvert_{\mbf{x}_J=\mbf{c}}))(\tau)=
\displaystyle\sum_{(x,y)\in \mal{A}_\tau(\mathbf{c})}	\zeta_q^{f_x-g_{y}},	
\end{equation}
where we define $\mal{A}_\tau(c)$ as
\begin{equation}\label{corr_restric_c}
\mal{A}_\tau(c)=\{(x,y)\in\mal{N}_c\times \mal{N}_c:0\leq x\leq L-\tau-1,y=x+\tau\}. 
\end{equation}
Below, we present another example to explain the above-described matters which we have mainly defined after Example \ref{exmo1}.\ccn
\begin{example}\label{ex2924}
Let us consider the same function $f$, as defined in (\ref{ankby1}), which is presented below:
$$f(x_{1},x_{2},x_{3},x_{4},x_{5})=2x_{1}x_{2}+4x_{2}x_{3}+x_{2}x_{4}+x_{2}x_{5}+3x_{1}x_{3}+2x_{4}x_{5}+x_{2}+2.$$ 
Again let us assume that $n_1=1$, $n_2=0$, and $J=\{J_1,J_2\}=\{j_{1,1}\}=\{2\}$, where $J_2=\emptyset$. Hence 
$\mbf{x}_J=({\mathbf{x}_1}_{J_1},{\mathbf{x}_2}_{J_2})=(x_{j_{1,1}})=(x_2)$. Then 
\begin{equation}\nonumber
\calV_{72}^0=\left\{(x_{1},x_{2},x_{3},x_{4},x_{5})\in \calV_{72}:x_{2}=0 \right\},
\end{equation}
and 
\begin{equation}\nonumber
\begin{split}
\mal{N}_0 & =\left\{x\in \mbf{Z}_{72}:x=9((x_{1},0,x_{3})\cdot (4,2,1))+(x_{4},x_{5})\cdot (3,1), x_{1}, x_{3}\in \mbf{Z}_2,~x_{4},x_{5}\in \mbf{Z}_3\right\}\\
&=\{x=2n,~2n+1:0\leq n\leq 35, ~n\!\!\!\!\mod\! 2=0\}.
%&=[0,17] \cup [36,53]. 
\end{split}
\end{equation} 
From (\ref{lbl1}), we have
\begin{equation}\nonumber
\begin{split}
\psi(f\arrowvert_{x_{2}=0})=&\left( \zeta_6^2     \zeta_6^2     \mbf{0}_2     \zeta_6^2     \zeta_6^5    \mbf{0}_2     \zeta_6^2     \zeta_6^2     \mbf{0}_2     \zeta_6^2     \zeta_6^5     \mbf{0}_2     \zeta_6^2     \zeta_6^2     \mbf{0}_2    \zeta_6^2     \zeta_6^5     \mbf{0}_2     \zeta_6^2     \zeta_6^2     \mbf{0}_2     \zeta_6^2     \zeta_6^5     \mbf{0}_2     \zeta_6^4     \zeta_6^4     \mbf{0}_2     \zeta_6^4     \zeta_6^1     \mbf{0}_2     \zeta_6^0     \zeta_6^0     \mbf{0}_2
\zeta_6^0     \zeta_6^3     \mbf{0}_2\right.\\& \left.    \zeta_6^2    \zeta_6^2     \mbf{0}_2     \zeta_6^2     \zeta_6^5     \mbf{0}_2     \zeta_6^0     \zeta_6^0     \mbf{0}_2     \zeta_6^0     \zeta_6^3     \mbf{0}_2    \zeta_6^4    \zeta_6^ 4     \mbf{0}_2     \zeta_6^4     \zeta_6^1     \mbf{0}_2\right),\end{split}
\end{equation}
where $\mbf 0_2=00$.
\ccn 
\end{example}
\ccn 
\subsection{Proposed $q$-ary Functions and Construction of CCCs for $k=2$}\label{contsec}
 In this subsection, we first provide all the details related to our proposed construction for $k=2$, and then, we shall extend the idea for $k\geq 3$. Let us assume that $f:\mbf Z_{p_1}^{m_1}\times \mbf Z_{p_2}^{m_2}\rightarrow \mbb Z_q$ be a $q$-ary function of Hamming degree $2$ over $m$-variables $x_1,x_2,\hdots,x_m$, where $m=m_1+m_2$, $m\geq 2$, $m_1,~m_2> 0$, and $q=p_1p_2$. Please note that $\mbf Z_{p_i}\subset \mbb Z_q$, and $\mbf Z_{p_1}^{m_1}\times \mbf Z_{p_2}^{m_2}\subset \mbb Z_q^m$, where $i=1,2$. For the one-to-one mappings $\pi:\{1,2,\hdots,m_1\}\rightarrow\{1,2,\hdots,m_1\}$ and $\pi':\{m_1+1,m_1+2,\hdots,m\}\rightarrow \{m_1+1,m_1+2,\hdots,m\}$, let us define $f$ as follows:
\begin{equation}\label{function10724}
	\begin{split}
f(\mbf x_1,\mbf x_2)&=f(x_1,x_2,\hdots,x_m)\\&=\frac{q}{p_1}\sum_{i=1}^{m_1-1} f_i(x_{\pi(i)})f_i'(x_{\pi(i+1)})+\sum_{\alpha=1}^{m_1}g_{\alpha}(x_\alpha)+\frac{q}{p_2}\sum_{j=1}^{m_2-1} h_j(x_{\pi'(m_1+j)}) h_j'(x_{\pi'(m_1+j+1)})\\&+\sum_{\beta=1}^{m_2}g'_{\beta} ({x_{\beta+m_1}})+\gamma f_0(x_{\pi(m_1)})h_0(x_{\pi'(m_1+1)}),
	\end{split}
\end{equation} 
where $f_i,~f_i'$, $h_j,~h_j',~f_0,~h_0,~g_{\alpha}$, and $g'_{\beta}$ are functions from $\mbb Z_q$ to $\mbb Z_q$. Now for $t=p_1t_2+t_1$, where $t_i\in \mbf Z_{p_i}$, let us define the following set of $q$-ary functions denoted by $C_t$ and the set of corresponding complex-valued sequences denoted by $\psi(C_t)$ as follows:\ccn
\begin{equation}\nonumber
	\begin{split}
	C_t=\left\{f+\frac{q}{p_1}(d_1x_{\pi(1)}+t_1x_{\pi(m_1)})+\frac{q}{p_2}(d_2x_{\pi'(m_1+1)}+t_2x_{\pi'(m)}):d_i\in\mbf Z_{p_i},~i=1,2  \right\},
	\end{split}
\end{equation}
and
\begin{equation}\label{ccset10724}
\begin{split}
\psi(C_t)=\left\{\psi\left(f+\frac{q}{p_1}(d_1x_{\pi(1)}+t_1x_{\pi(m_1)})+\frac{q}{p_2}(d_2x_{\pi'(m_1+1)}+t_2x_{\pi'(m)})\right):d_i\in\mbf Z_{p_i},~i=1,2  \right\}.
\end{split}
\end{equation}
\begin{theorem}\label{th314824}
	Let $\mal C=\{\psi(C_t):t=0,1,\hdots,p_1p_2-1\}$ be the set of complex-valued codes defined in (\ref{ccset10724}). Then for any choices of $f_0,~h_0,~g_{\alpha},~g'_{\beta}$, one-to-one mappings $\pi,~\pi'$, and $\lambda\in \mbb Z_q$, the code set $\mal C$ forms $(p_1p_2,L)$-CCCs iff $f_i$ and $f_i'$ permute the set $\mbf Z_{p_1}$ under modulo $p_1$, and $h_j$ and $h_j'$ permute the set $\mbf Z_{p_2}$ under modulo $p_2$,\ccn where the index terms $i,~j,~\alpha$, and $\beta$ follow the same defintion as defined in (\ref{function10724}), and $L=p_1^{m_1}p_2^{m_2}$.	
\end{theorem}
\begin{IEEEproof}
For $t=p_1t_2+t_1$ and $t'=p_1t_2'+t_1'$, let $F_{x,t}=f_x+\frac{q}{p_1}(t_1x_{\pi(m_1)})+\frac{q}{p_2}(t_2x_{\pi'(m)})$, and $F_{x,t'}=f_x+\frac{q}{p_1}(t_1'x_{\pi(m_1)})+\frac{q}{p_2}(t_2'x_{\pi'(m)})$. With this set of notations, for $0<\tau<L$, the ACCF between $\psi(C_t)$ and $\psi(C_{t'})$ can be expressed as 
\begin{equation}\label{corr1_10724}
	\begin{split}
	\Theta(\psi(C_t),\psi(C_{t'}))(\tau)=\sum_{d_1,d_2}\Theta\left(\psi\left(F_{x,t}+\frac{q}{p_1}(d_1x_{\pi(1)})+\frac{q}{p_2}(d_2x_{\pi'(m_1+1)})\right),\right.\\ \left. \psi\left(F_{x,t'}+\frac{q}{p_1}(d_1x_{\pi(1)})+\frac{q}{p_2}(d_2x_{\pi'(m_1+1)})\right)\right)(\tau).
	\end{split}
\end{equation}
For $0<\tau<L$, $0\leq x<L-\tau$ and $y=x+\tau$, let us assume $(x_1,x_2,\hdots,x_{m})$ and $(y_1,y_2,\hdots,y_m)$ are their vector representations, respectively. Then
\begin{equation}\label{corr2_10724}
	\begin{split}
\Theta &\left(\psi\left(F_{x,t}+\frac{q}{p_1}(d_1x_{\pi(1)})+\frac{q}{p_2}(d_2x_{\pi'(m_1+1)})\right), \psi\left(F_{x,t'}+\frac{q}{p_1}(d_1x_{\pi(1)})+\frac{q}{p_2}(d_2x_{\pi'(m_1+1)})\right)\right)(\tau)\\
&= \sum_{x=0}^{L-\tau-1} \xi_q^{F_{x,t}-F_{y,t'}} \xi_{p_1}^{d_1(x_{\pi(1)}-y_{\pi(1)})} \xi_{p_2}^{d_2(x_{\pi'(m_1+1)}-y_{\pi'(m_1+1)})},
	\end{split}
\end{equation}
where $F_{x,t}=f_x+\frac{q}{p_1}(t_1x_{\pi(m_1)})+\frac{q}{p_2}(t_2x_{\pi'(m)})$ and $F_{y,t'}=f_y+\frac{q}{p_1}(t_1'y_{\pi(m_1)})+\frac{q}{p_2}(t_2'y_{\pi'(m)})$.
Applying (\ref{corr2_10724}) in (\ref{corr1_10724}), we have
\begin{equation}\label{corr3_10724}
	\begin{split}
	\Theta(\psi(C_t),\psi(C_{t'}))(\tau)= \sum_{x=0}^{L-\tau-1} \xi_q^{F_{x,t}-F_{y,t'}}\left(\sum_{d_1} \xi_{p_1}^{d_1(x_{\pi(1)}-y_{\pi(1)})} \right) \left(\sum_{d_2}\xi_{p_2}^{d_2(x_{\pi'(m_1+1)}-y_{\pi'(m_1+1)})}\right).
	\end{split}
\end{equation}
To derive (\ref{corr3_10724}), we can go through the following four cases:
\begin{enumerate}
	\item $x_{\pi(1)}=y_{\pi(1)},~x_{\pi'(m_1+1)}=y_{\pi'(m_1+1)}$.
	\item $x_{\pi(1)}\neq y_{\pi(1)},~x_{\pi'(m_1+1)}=y_{\pi'(m_1+1)}$.
	\item $x_{\pi(1)}= y_{\pi(1)},~x_{\pi'(m_1+1)}\neq y_{\pi'(m_1+1)}$.
	\item $x_{\pi(1)}\neq y_{\pi(1)},~x_{\pi'(m_1+1)}\neq y_{\pi'(m_1+1)}$.
	\end{enumerate}
Among the above four cases, for the last three cases, from (\ref{corr3_10724}), we have $\Theta(\psi(C_t),\psi(C_{t'}))(\tau)=0$. Then at this stage, it says, we only need to take care of the first case $x_{\pi(1)}=y_{\pi(1)},~x_{\pi'(m_1+1)}=y_{\pi'(m_1+1)}$.
In this case (\ref{corr3_10724}) reduces to the following:
\begin{equation}\label{corr4_10724}
\begin{split}
\Theta(\psi(C_t),\psi(C_{t'}))(\tau)= p_1p_2\sum_{x=0}^{L-\tau-1} \xi_q^{F_{x,t}-F_{y,t'}}.
\end{split}
\end{equation}	
Let $1<v\leq m_1$ is the smallest positive integer such that $x_{\pi(v)}\neq y_{\pi(v)}$.
Again let us define another integer $x^u$ and its vector representation is $(x_1,x_2,\hdots,(x_{\pi(v-1)}+u)\mod q,\hdots,x_{m_1},\hdots,x_m)$, where $1\leq u<p_1$. It says that the vector representations of $x$ and $x^u$ differ only at the $\pi(v-1)$th position. Similarly, we define another integer $y^u$ having vector representation $(y_1,y_2,\hdots,(y_{\pi(v-1)}+u)\mod q,\hdots,y_{m_1},\hdots,y_m)$ which differ only at the position $\pi(v-1)$ with the vector representation of $y$. Then
\begin{equation}\nonumber
	\begin{split}
	F&_{x^u,t}-F_{x,t}\\=&\frac{q}{p_1}\left[f_{v-2}(x_{\pi(v-2)})f'_{v-2}((x_{\pi(v-1)}+u)\!\!\!\!\!\mod q) +f_{v-1}((x_{\pi(v-1)}+u)\!\!\!\!\!\mod q)f'_{v-1}(x_{\pi(v)})\right.\\& \left.-       f_{v-2}(x_{\pi(v-2)})f'_{v-2}(x_{\pi(v-1)}) - f_{v-1}(x_{\pi(v-1)})f'_{v-1}(x_{\pi(v)})\right]\\
                & +g_{\pi(v-1)}((x_{\pi(v-1)}+u)\!\!\!\!\!\mod q)-g_{\pi(v-1)}(x_{\pi(v-1)}),
	\end{split}
\end{equation}
and

\begin{equation}\nonumber
\begin{split}
F&_{y^u,t'}-F_{y,t'}\\=&\frac{q}{p_1}\left[f_{v-2}(x_{\pi(v-2)})f'_{v-2}((x_{\pi(v-1)}+u)\!\!\!\!\!\mod q) +f_{v-1}((x_{\pi(v-1)}+u)\!\!\!\!\!\mod q)f'_{v-1}(y_{\pi(v)})\right.\\& \left.-       f_{v-2}(x_{\pi(v-2)})f'_{v-2}(x_{\pi(v-1)}) - f_{v-1}(x_{\pi(v-1)})f'_{v-1}(y_{\pi(v)})\right]\\
& +g_{\pi(v-1)}((x_{\pi(v-1)}+u)\!\!\!\!\!\mod q)-g_{\pi(v-1)}(x_{\pi(v-1)}).
\end{split}
\end{equation}  	
Now 
\begin{equation}\nonumber
\begin{split}
F&_{x^u,t}-F_{y^u,t'}-(F_{x,t}-F_{y,t'})\\&=\frac{q}{p_1}(f'_{v-1}(x_{\pi(v)})-f'_{v-1}(y_{\pi(v)}))(f_{v-1}((x_{\pi(v-1)}+u)\!\!\!\!\!\mod q)-f_{v-1}(x_{\pi(v-1)})).
\end{split}
\end{equation}	
Then 
\begin{equation}\label{cond11724}
	\begin{split}
	\sum_{u=1}^{p_1-1} \xi_q^{F_{x^u,t}-F_{y^u,t'}-(F_{x,t}-F_{y,t'})}=-1+\sum_{u=0}^{p_1-1} \xi_{p_1}^{rp(u)},
	\end{split}
\end{equation}
where $r=f'_{v-1}(x_{\pi(v)})-f'_{v-1}(y_{\pi(v)})$ and $p(u)=f_{v-1}((x_{\pi(v-1)}+u)\!\!\!\mod q)-f_{v-1}(x_{\pi(v-1)})$. 
\noindent \textit{\textbf{Sufficiency:}}
Let $f_i$ and $f_i'$ permute $\mbf Z_{p_1}$ under modulo $p_1$, where $i=1,2,\hdots,m_1-1$. \ccn
%Then $r\neq 0\mod p_1$, and as $p_1|q$, it implies $r\neq 0\mod p_1$. For $u_1\neq %u_2\in\mbf Z_{p_1}$, it is true that $p(u_1)\neq p(u_2)\mod q$. As $p_1|q$,  $p(u_1)\neq %p(u_2)\mod q$ implies $p(u_1)\neq p(u_2)\mod p_1$. 
Therefore, $\sum_{u=0}^{p_1-1} \xi_{p_1}^{rp(u)}=0$. 
Then form (\ref{cond11724}), we have
\begin{equation}\label{cond2_11724}
\begin{split}
\xi_{p_1}^{F_{x,t}-F_{y,t'}}+\sum_{u=1}^{p_1-1} \xi_q^{F_{x^u,t}-F_{y^u,t'}}=0.
\end{split}
\end{equation} 
Now combining all the above results, it is clear that 
\begin{equation}\nonumber
	\Theta(\psi(C_t),\psi(C_{t'}))(\tau)=0~\text{if}~f_i~ \text{and}~ f_i'~ \text{permute}~\mbf Z_{p_1}~\textnormal{under modulo }~p_1,~i=1,2,\hdots,m_1-1.
\end{equation}\ccn
%$\Theta(\psi(C_t),\psi(C_{t'}))(\tau)=0$ iff $f_i$ and $f_i'$ permute $\mbf Z_{p_1}$, $i=1,2,\hdots,m_1-1$.	
When $m_1+1<v\leq m$, using similar argument as drawn in this proof, we can reach the following result:
\begin{equation}\nonumber
\Theta(\psi(C_t),\psi(C_{t'}))(\tau)=0~\text{if}~h_j~ \text{and}~ h_j'~ \text{permute}~\mbf Z_{p_2}~\textnormal{under modulo}~p_2,~j=1,2,\hdots,m_2-1.
\end{equation}\ccn
%\end{case}
Again combining all the above results, we have
\begin{equation}\label{cond3_11724}
\begin{split}
\Theta(\psi(C_t),\psi(C_{t'}))(\tau)=0~\text{if}~&f_i,~ f_i'~ \text{permute}~\mbf Z_{p_1}~\textnormal{under modulo }~p_1,~\text{and}~h_j,~ h_j'~ \text{permute}~\mbf Z_{p_2}\\&\textnormal{under modulo }~p_2,
\end{split}
\end{equation}
where $i=1,2,\hdots,m_1-1$ and $j=1,2,\hdots,m_2-1$. Similarly, we can reach to the same result as in (\ref{cond3_11724}) for $\tau<0$. 
For $\tau=0$, from (\ref{corr4_10724}), we have
\begin{equation}
	\begin{split}
	\Theta(\psi(C_t),\psi(C_{t'}))(0)=p_1p_2\sum_{x=0}^{L-1} \xi_{p_1}^{(t_1-t_1')x_{\pi(m_1)}} \xi_{p_2}^{(t_2-t_2')x_{\pi(m)}}=
	\begin{cases}
          p_1^{m_1+1}p_2^{m_2+1},& t=t',\\
          0,& \text{otherwise}.
	\end{cases}
	\end{split}
\end{equation}
Again combining all the above results, it is clear the set of codes $\mal C$ forms $(p_1p_2,L)$-CCCs if $f_i,~f_i'$ permute the set $\mbf Z_{p_1}$~\textnormal{under modulo }~$p_1$, and $h_j,~h_j'$ permute the set $\mbf Z_{p_2}$~\textnormal{under modulo }~\(p_2\).\ccn

\noindent \textit{\textbf{Necessity:}}
Let us assume that the code set $C$ forms $(p_1p_2,p_1^{m_1}p_2^{m_2})$-CCCs. Here, our goal is to show that the functions $f_i,~f_i':\mbb Z_{q}\rightarrow \mbb Z_{q}$ permute  the set $\mbf Z_{p_1}$ under modulo $p_1$, and the functions $h_j,~h_j':\mbb Z_{q}\rightarrow \mbb Z_{q}$ permute $\mbf Z_{p_2}$ under modulo $p_2$, where $i=1,2,\hdots,m_1$ and $j=1,2,\hdots,m_2$. Without loss of generality, let us assume that $\pi$ are $\pi'$ are the identity permutations and $\tau=p_1^{m_1}-np_1^{m_1-i}$, where $1\leq n<p_1$. Let $(x_1,x_2,\hdots,x_m)$ is the $q$-ary vector representation of $x$, i.e.,  $x=\sum_{j_1=1}^{m_1} x_j p_1^{j_1-1}+p_1^{m_1}\sum_{j_2=1}^{m_2}x_{m_1+j}p_2^{j_2-1}$. Then $y=x+\tau$, implies 
\begin{equation}\label{argu112824}
\begin{split}
y=x+\tau=&\sum_{j_1=1}^{m_1-i} x_j p_{1}^{j_1-1}+(x_{m_1-i+1}+p_1-n)p_1^{m-i}+\sum_{j_1=m_1-i+2}^{m_1} (x_j+p_1-1) p_1^{j_1-1}\\&+p_1^{m_1}\sum_{j_2=1}^{m_2}x_{m_2+j}p_2^{j_2-1}.
\end{split}
\end{equation}
For $(y_1,y_2,\hdots,y_m)$ as the vector representation of $y$, 
from (\ref{argu112824}), we have $x_{j_1}=y_{j_1}$, for $j_1=1,2,\hdots,m-i$, $y_{m-i+1}=x_{m-i+1}+q-n$, where $x_{m-i+1}=0,1,\hdots,n-1$. Besides, $y_{j_1'}=x_{j_1'}+q-1$, where $x_{j_1'}=0$, for $j_1'=m-i+2,m-i+3,\hdots,m-1$, and for $j_2=1,2,\hdots,m_2$, $x_{j_2}=y_{j_2}$. Let $S_n$ be the set of all the pairs $(x,y)$ with $y-x=p_1^{m_1}-np_1^{m_1-i}$, can be expressed as 
\begin{equation}\label{argu212824}
\begin{split}
S_n=&\left\{(x,y): y=x+p_1^{m_1}-np_1^{m_1-i},~ y_{j_1}=x_{j_1},~y_{m_1-i+1}=x_{m_1-i+1}+p_1-n,~ \right. \\& \left.y_{j_1'}=x_{j_1'}+p_1-1, x_{j_1}=0,1,\hdots,p_1-1, x_{m_1-i+1}=0,1,\hdots,n-1, ~x_{j_1'}=0, \right.\\&\left. \text{and},~y_{m_1+j_2}=x_{m_1+j_2}, ~\text{where}~ j_1=1,2,\hdots,m_1-i,~j_1'=m-i+2,\hdots,m-1,\right.\\&\left.\text{and}~j_2=1,2,\hdots,m_2\right\}. 
\end{split}
\end{equation}
Then $|S_n|=np_2^{m_2}p_1^{m_1-i}$, and for $1\leq n<p_1$, (\ref{corr4_10724}) can be represented as
\begin{equation}\label{key112824}
\begin{split}
\Theta(\psi(C_t),\psi(C_{t'}))(q^m-nq^{m-i})&=\sum_{(x,y)\in S_n} \xi_q^{F_{x,t}-F_{y,t'}}\\
%&=\sum_{(x,y)\in \cup_{\alpha=1}^{q^{m-i-1}} S_{1,\alpha}} \xi_q^{F_{x,t}-F_{y,t'}}\\
&=p_1p_2\sum_{\beta =1}^{p_2^{m_2} }\sum_{\alpha=1}^{p_1^{m_1-i-1}}\sum_{x_{\alpha,m_1-i+1}=0}^{n-1} \sum_{u=0}^{p_1-1} \xi_q^{F_{x_{\alpha,\beta}^u,t}-F_{y_{\alpha,\beta}^u,t'}}, 
\end{split}
\end{equation}
where 
$$x_{\alpha,\beta}^u=(x_{\alpha,1},x_{\alpha,2},\hdots,x_{\alpha,m_1-i-1},u,x_{\alpha,m_1-i+1},0,\hdots,0,x_{\beta,m_1+1},\hdots,x_{\beta,m}),$$ and 
$$y_{\alpha,\beta}^u=(x_{\alpha,1},x_{\alpha,2},\hdots,x_{\alpha,m_1-i-1},u,x_{\alpha,m-i+1}+q-n,q-1,\hdots,q-1,x_{\beta,m_1+1},\hdots,x_{\beta,m}).$$
%Now $F_{x_{\alpha,\beta}^u,t}$
Then
\begin{equation}\label{aug113824}
	\begin{split}
&F_{x_{\alpha,\beta}^u,t}-F_{y_{\alpha,\beta}^u,t'}\\
=&\frac{q}{p_1}\left[f_{m-i}(u)\left(f_{m-i}'(x_{\alpha,m_1-i+1})-f_{m-i}'(x_{\alpha,m_1-i+1}+p_1-n)\right)+f_{m-i+1}(x_{\alpha,m-i+1})f_{m-i+1}'(0)\right.\\&\left.-f_{m-i+1}(x_{\alpha,m-i+1}+p_1-n)f_{m-i+1}'(p_1-1)\right]\!+\!\frac{q}{p_1}\!\!\sum_{j=m_1-i+2}^{m_1}\!\!\!\!\!\!\! \left(f_j(0)f_j'(0)-f_j(p_1-1)f_j'(p_1-1)\right)\\&+\left(g_{m_1-i+1}(x_{\alpha,m_1-i+1})-g_{m_1-i+1}(x_{\alpha,m_1-i+1}+p_1-n)\right)+
\sum_{j=m_1-i+2}^{m_1}\left(g_j(0)-g_j(p_1-1)\right)\\&+\lambda h_0(x_{\beta,m_1+1})\left(f_0(0)-f_0(p_1-1)\right)+\frac{q}{p_1}(t_1\cdot 0-t_1'\cdot (p_1-1))+\frac{q}{p_2}x_{\beta,m}(t_2-t_2'),
	\end{split}
\end{equation}
and
\begin{equation}\label{aug213824}
\begin{split}
F&_{x_{\alpha,\beta}^u,t}-F_{y_{\alpha,\beta}^u,t'}-(F_{x_{\alpha,\beta}^0,t}-F_{y_{\alpha,\beta}^0,t'})\\&=\frac{q}{p_1}\left(f_{m-i}'(x_{\alpha,m_1-i+1})-f_{m-i}'(x_{\alpha,m_1-i+1}+p_1-n)\right)\left(f_{m-i}(u)-f_{m-i}(0)\right).
\end{split}
\end{equation}
Now let us assume that $n=1$, which implies $x_{\alpha,m_1-i+1}=0$.  Using  (\ref{aug113824})
and (\ref{aug213824}), (\ref{key112824}) can be represented as 
\begin{equation}\label{aug313824}
\begin{split}
\Theta&(\psi(C_t),\psi(C_{t'}))(q^m-nq^{m-i})\\&=\sum_{(x,y)\in S_n} \xi_q^{F_{x,t}-F_{y,t'}}\\
%&=\sum_{(x,y)\in \cup_{\alpha=1}^{q^{m-i-1}} S_{1,\alpha}} \xi_q^{F_{x,t}-F_{y,t'}}\\
&=p_1p_2^{m_2-1}\sum_{x_{\beta,m} =0}^{p_2-1 } \sum_{x_{\beta,m_1+1}=0}^{p_2-1}\sum_{\alpha=1}^{p_1^{m_1-i-1}} \xi_q^{F_{x_{\alpha,\beta}^0,t}-F_{y_{\alpha,\beta}^0,t'}} \left(\sum_{u=0}^{p_1-1} \xi_{p_1}^{r\left(f_{m-i}(u)-f_{m-i}(0)\right)}\right), 
\end{split}
\end{equation}
where $r=f_{m-i}'(0)-f_{m-i}'(p_1-1)$. As $n=1$ and $x_{\alpha,m_1-i+1}=0$, from (\ref{aug113824}), it can be observed that
$F_{x_{\alpha_1,\beta}^u,t}-F_{y_{\alpha_1,\beta}^u,t'}=F_{x_{\alpha_2,\beta}^u,t}-F_{y_{\alpha_2,\beta}^u,t'}=F_{x_{\alpha,\beta}^u,t}-F_{y_{\alpha,\beta}^u,t'}$ for all $0\leq \alpha_1,~\alpha_2\leq p_1^{m_1-i-1}$. Then (\ref{aug313824}) can be represented as 
\begin{equation}\label{aug413824}
\begin{split}
\Theta&(\psi(C_t),\psi(C_{t'}))(q^m-nq^{m-i})\\&=\sum_{(x,y)\in S_n} \xi_q^{F_{x,t}-F_{y,t'}}\\
%&=\sum_{(x,y)\in \cup_{\alpha=1}^{q^{m-i-1}} S_{1,\alpha}} \xi_q^{F_{x,t}-F_{y,t'}}\\
&=p_1^{m_1-i} p_2^{m_2-1}\left(\sum_{x_{\beta,m} =0}^{p_2-1 } \sum_{x_{\beta,m_1+1}=0}^{p_2-1} \xi_q^{F_{x_{\alpha,\beta}^0,t}-F_{y_{\alpha,\beta}^0,t'}}\right) \left(\sum_{u=0}^{p_1-1} \xi_{p_1}^{r\left(f_{m-i}(u)-f_{m-i}(0)\right)}\right).
\end{split}
\end{equation}
Again for $n=1$, (\ref{aug113824}) can be expressed as 
\begin{equation}\label{aug513824}
\begin{split}
F_{x_{\alpha,\beta}^u,t}-F_{y_{\alpha,\beta}^u,t'}=R+
\lambda h_0(x_{\beta,m_1+1})\left(f_0(0)-f_0(p_1-1)\right)+\frac{q}{p_2}x_{\beta,m}(t_2-t_2'),
\end{split}
\end{equation}
where
\begin{equation}\nonumber
	\begin{split}
R&=\frac{q}{p_1}\left(f_{m-i}(u)\left(f_{m-i}'(0)-f_{m-i}'(p_1-1)\right)\right)+\frac{q}{p_1}\sum_{j=m_1-i+1}^{m_1} \left(f_j(0)f_j'(0)-f_j(p_1-1)f_j'(p_1-1)\right)\\&+
\sum_{j=m_1-i+1}^{m_1}\left(g_j(0)-g_j(p_1-1)\right)+\frac{q}{p_1}(t_1\cdot 0-t_1'\cdot (p_1-1)).
	\end{split}
\end{equation}
Then (\ref{aug413824}) can be expressed as
\begin{equation}\label{aug613824}
\begin{split}
\Theta&(\psi(C_t),\psi(C_{t'}))(q^m-nq^{m-i})\\&=\sum_{(x,y)\in S_n} \xi_q^{F_{x,t}-F_{y,t'}}\\
&=\xi_q^R p_1^{m_1-i} p_2^{m_2-1}\left(\sum_{x_{\beta,m_1+1} =0}^{p_2-1 } \xi_q^{\lambda h_0(x_{\beta,m_1+1})(f_0(0)-f_0(p_1-1))}\right)\left( \sum_{x_{\beta,m}=0}^{p_2-1} \xi_{p_2}^{x_{\beta,m}(t_2-t_2')}\right) \\&~~~~~~~~~~~~~~~~~~~~~\left(\sum_{u=0}^{p_1-1} \xi_{p_1}^{r\left(f_{m-i}(u)-f_{m-i}(0)\right)}\right).
\end{split}
\end{equation}
Now our assumption is that for any choices of $f_0,~g_{\alpha'}:\mbb Z_{q}\rightarrow \mbb Z_q$, and $h_0, ~g'_{\beta'}:\mbb Z_{q}\rightarrow \mbb Z_q$, $\mal C$ forms $(p_1p_2,L)$-CCCs, where $\alpha'=1,2\hdots,m_1$ and $\beta'=1,2,\hdots,m_2$. Then without loss of generality, we again assume $t_2=t_2'$ and $f_0(0)=f_0(p_1-1)$. With this assumption, from (\ref{aug613824}), we have
\begin{equation}\nonumber
	\begin{split}
\Theta(\psi(C_t),\psi(C_{t'}))(q^m-nq^{m-i})=\xi_q^R p_1^{m_1-i} p_2^{m_2+1}\sum_{u=0}^{p_1-1}\xi_{p_1}^{r\left(f_{m-i}(u)-f_{m-i}(0)\right)}=0,	
	\end{split}
\end{equation}
 which implies $r\neq 0 \mod p_1$ and $f_{m-i}:\mbb Z_{q}\rightarrow \mbb Z_{q} $ permutes the set $\mbf Z_{p_1}$ under the operation modulo $p_1$, where $i=1,2,\hdots, m_1-1$. Now $r\neq 0$ implies $f'_{m-i}(0)\neq f'_{m-i}(p_1-1)$, and using similar approach as in the \textit{\textbf{Necessity}} part in the proof of \textbf{Theorem \ref{th824}}, we can show that $f'_{m-i}(u_1)\neq f'_{m-i}(u_2)\mod p_1$ for all $u_1\neq u_2\in \mbf Z_{p_1}$, which implies $f'_{m-i}:\mbb Z_{q}\rightarrow \mbb Z_{q}$ also permutes 
$\mbf Z_{p_1}$ under modulo $p_1$. Similarly as above, it can also be shown that the functions $h_j,~h_j':\mbb{Z}_{q}\rightarrow \mbb Z_{q}$ permute the set $\mbf Z_{p_2}$ under the operation modulo $p_2$, where $j=1,2,\hdots,m_2-1$.
\end{IEEEproof}
In \textbf{Theorem \ref{th824}}, we have proposed $(q,q^m)$-CCCs, and then in \textbf{Corollary \ref{cor15824}}, we have extended the result $(q,q^m)$-CCCs to $(q^{n+1},q^m)$-CCCs with the idea of restricting a function over $n$-variables. Similarly, with the help of \textbf{Theorem \ref{th314824}} and applying the idea of restricting a function over some variables, we shall extend the result $(p_1p_2,p_1^{m_1}p_2^{m_2})$-CCCs to $(p_1^{n_1+1}p_2^{n_2+1},p_1^{m_1}p_2^{m_2})$-CCCs.

Let us assume that after restricting $f$ at $\mbf x_J=\mbf c$, it results in a $q$-ary function $f\arrowvert_{\mathbf{x}_J=\mathbf{c}}$ of Hamming degree $2$ over rest of the $m-n$ variables labelled, $\{1,2,\hdots,m\}\setminus J_1\cup J_2$, where we recall the notations for $k=2$, defined in Section \ref{sectn5} as  $\mbf x=(x_1,x_2,\hdots,x_m)=(\mbf x_1,\mbf x_2)$, $\mbf x_J=({\mbf x_1}_{J_1},{\mbf x_2}_{J_2})$, and $\mbf c=(\mbf c_1,\mbf c_2)\in \mbf Z_{p_1}^{n_1}\times \mbf Z_{p_1}^{n_2}$, where $n_1+n_2=n$, $0\leq n_1\leq m_1-1$, and $0\leq n_2\leq m_2-1$, $|J|=|\{J_1,J_2\}|=n_1+n_2=n$, $\mbf x_1=(x_1,x_2,\hdots,x_{m_1})$ and 
$\mbf x_2=(x_{m_1+1},x_{m_1+2},\hdots,x_m)$. As in this section, we deal for $k=2$, let us re-define $J_1$ and $J_2$ with less complex notations as   
$J_1=(j_{1},j_{2},\hdots,j_{n_1})\subset \{1,2,\hdots,m_1\}$ and $J_2=(j_1',j_2',\hdots,j_{n_2}')\subset \{m_1+1,m_1+2,\hdots,m\}$. Let  $f\arrowvert_{\mathbf{x}_J=\mathbf{c}}$ is in the following form:
\begin{equation}\label{fun22824}
	\begin{split}
	f\arrowvert_{\mathbf{x}_J=\mathbf{c}}&=\frac{q}{p_1}\sum_{i=1}^{m_1-n_1-1} f_i(x_{\pi_c(i)})f_i'(x_{\pi_c(i+1)})+\sum_{\alpha=1}^{m_1-n_1} g_\alpha(x_{\pi_c(\alpha)})\\&
+\frac{q}{p_2}	\sum_{j=1}^{m_2-n_2-1} h_j(x_{\pi_c'(j)})h_{j}'(x_{\pi_c'(j+1)})+\sum_{\beta=1}^{m_2-n_2} g'_{\beta}(x_{\pi_c'(\beta)})+\gamma f(x_{\pi_c(m_1-n_1)})h(x_{\pi_c'(1)}),
	\end{split}
\end{equation}
where $f_i,~f_i',~h_j,~h_j',~f, g_{\alpha},~h$, and $~g_{\beta}'$ are functions from $\mbb Z_q$ to $\mbb Z_q$, and $\lambda\in \mbb Z_q$. For $c=c_1\Delta_1'+c_2\Delta_2'\in \mbf{Z}_{L'}$, where $c_i=\sum_{j=1}^{n_i} c_{i,j} p_i^{j-1}\in \mbf Z_{L_i'}$, and $\Delta_i'=\frac{L'}{L_i'L_{i+1}'\cdots L_k'}$, $L'=p_1^{n_1}p_2^{n_2}$, $L_i'=p_i^{n_i}$, $i=1,2$, $\pi_c:\{1,2,\hdots,m_1-n_1\}\rightarrow \{1,2,\hdots,m_1\}\setminus J_1$ and $\pi_c':\{1,2,\hdots,m_2-n_2\}\rightarrow \{m_1+1,m_1+2,\hdots,m\}\setminus J_2$ are one-to-one mappings. For $t=t_2 p_1^{n_1+1}+t_1$, where $t_i\in \mbf Z_{p_i^{n_i+1}}$, let us assume $(t_{1,1},t_{1,2},\hdots,t_{1,n_1+1})\in \mbf Z_{p_1}^{n_1+1}$ and $(t_{2,1},t_{2,2},\hdots,t_{2,n_2+1})\in\mbf Z_{p_2}^{n_2+1}$ are the vector representation of $t_1$ and $t_2$, respectively. Let us define the following set $C_t$ containing $p_1^{n_1+1}p_2^{n_2+1}$ $q$-ary functions as
\begin{equation}\label{key114824}
	\begin{split}
	C_t=&\left\{f+ \sum_{i=1}^2\frac{q}{p_i}\left(\mbf d_i+\mbf t_i\right)\cdot {\mbf x_i }_{J_i} + \frac{q}{p_1}\left(d_{1,n_1+1}x_{\pi(1)} +t_{1,n_1+1}x_{\pi(m_1-n_1)} \right)
\right.	\\& \left.+\frac{q}{p_2}\left(d_{2,n_2+1}x_{\pi'(1)} +t_{2,n_2+1}x_{\pi'(m_2-n_2)} \right): \mbf d_i\in \mbf Z_{p_i}^{n_i}, d_{i,n_i+1}\in \mbf Z_{p_i}, i=1,2 \right\},
	\end{split}
\end{equation}
where $\mbf d_i=(d_{i,1},d_{i,2},\hdots,d_{i,n_i})$, $\mbf t_i=(t_{i,1},t_{i,2},\hdots,t_{i,n_i})$. Besides 
$\pi=\pi_c$ and $\pi'=\pi_c'$ if
 $\mbf x_J=\mbf c$. 
\begin{corollary}\label{corr160824}
Let $\mal C=\{\psi(C_t):0\leq t<p_1^{n_1+1}p_2^{n_2+1} \}$ be the set of codes corresponding to the set of $q$-ary functions defined in (\ref{key114824}). Then for any choice of $f,~h,~g_\alpha,~g_\beta'$, one-to-one mappings $\pi_c,~\pi_c'$ and $\lambda\in \mbb Z_q$ in (\ref{fun22824}), the code set $\mal C$ forms 
$(p_1^{n_1+1}p_2^{n_2+1}, L)$-CCCs iff $f_i,~f_i'$ permute the set $\mbf Z_{p_1}$ under modulo $p_1$, and $h_j,~h_j'$ permute the set $\mbf Z_{p_2}$ under modulo $p_2$, where the index terms $i,~j,~\alpha$, and $\beta$
follow the same definition as defined in (\ref{key114824}), and $L=p_1^{m_1}p_2^{m_2}$. 
\end{corollary}
\begin{IEEEproof}
For any two integers $t=t_2p_1^{n_1+1}+t_1$ and $t'=t_2' p_1^{n_1+1}+t_1'$, let $\mbf t_1\in \mbf Z_{p_1^{n_1+1}}$ and $\mbf t_2\in \mbf Z_{p_2^{n_2+1}}$ are the vector representations of $t_1$ and $t_2$, and $\mbf t_1'\in Z_{p_1^{n_1+1}}$ and $\mbf t_2'\in \mbf Z_{p_2^{n_2+1}}$ are the vector representations of $t_1'$ and $t_2'$, respectively.
We also assume
$F_{x,t}=f+ \frac{q}{p_1} (t_{1,n_1+1}x_{\pi(m_1-n_1)})+\frac{q}{p_2}(t_{2,n_2+1}x_{\pi'(m_2-n_2)})$
and $F_{x,t'}=f+ \frac{q}{p_1} (t_{1,n_1+1}'x_{\pi(m_1-n_1)})+\frac{q}{p_2}(t_{2,n_2+1}' x_{\pi'(m_2-n_2)})$. For $0<\tau<L$, 
the ACCF between $\psi(C_t)$ and $\psi(C_{t'})$ can be expressed as 
\begin{equation}\label{key214824}
	\begin{split}
	\Theta&(\psi(C_t),\psi(C_{t'}))(\tau)\\=&p_1^{n_1}p_2^{n_2}\sum_{\mbf c\in \mal V_{L'}}
	\xi_{p_1}^{(\mbf t_1-\mbf t_1')\cdot c_1} \xi_{p_2}^{(\mbf t_2-\mbf t_2')\cdot c_2} \sum_{\substack{d_{1,n_1+1}\\d_{2,n_2+1}} } \Theta\left(\psi\left( F_{x,t}+\frac{q}{p_1}(d_{1,n_1+1}x_{\pi(1)})+\frac{q}{p_2}(d_{2,n_2+1}x_{\pi'(1)})\right)\arrowvert_{\mathbf{x}_J=\mathbf{c}},\right. \\&\left. \psi\left(F_{x,t}+\frac{q}{p_1}(d_{1,n_1+1}x_{\pi(1)})+\frac{q}{p_2}(d_{2,n_2+1}x_{\pi'(1)})\right)\arrowvert_{\mathbf{x}_J=\mathbf{c}} \right)(\tau).
	\end{split}
\end{equation}
It can be observed from (\ref{key214824}) that, similarly as the proof of \textbf{Theorem \ref{th314824}}, it can be shown that the code set $\mal C$ forms $(p_1^{n_1+1}p_2^{n_2+1},L)$-CCCs iff $f_i,~f_i'$ permute the set $\mbf Z_{p_1}$ under modulo $p_1$, and $h_j,~h_j'$ permute the set $\mbf Z_{p_2}$ under modulo $p_2$.
\end{IEEEproof}
\begin{example}
		\begin{figure}[t]
		\includegraphics[width=10cm]{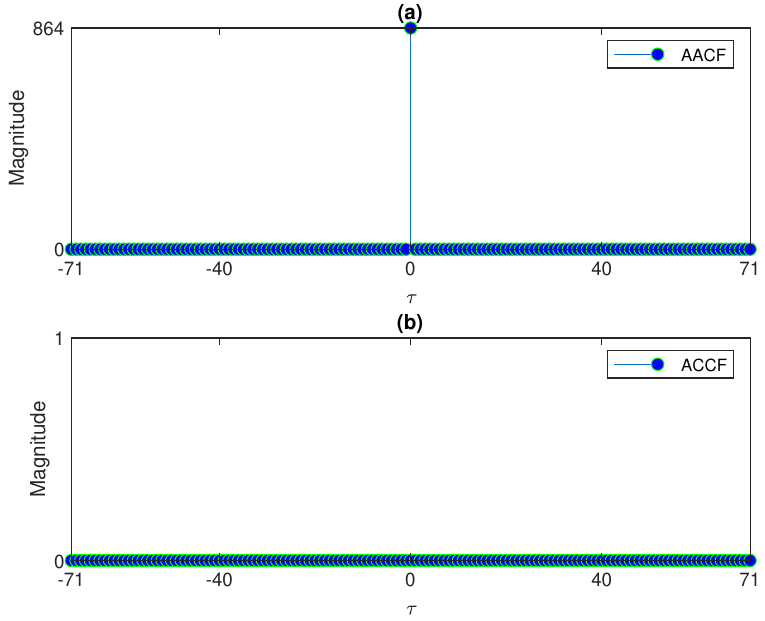}
		\caption{Correlation Plotting for the codes $\psi(C_{1})$ and $\psi(C_{11})$.}
		\label{ankut111}
	\end{figure}
Let us consider the same function as in Example \ref{ex2924}. Then 
$f(x_{1},x_{2},x_{3},x_{4},x_{5})=2x_{1}x_{2}+4x_{2}x_{3}+x_{2}x_{4}+x_{2}x_{5}+3x_{1}x_{3}+2x_{4}x_{5}+x_{2}+2$, 
$n_1=1$, $n_2=0$, and $J=\{J_1,J_2\}=\{j_{1,1}\}=\{2\}$, where $J_2=\emptyset$. Hence 
$\mbf{x}_J=({\mathbf{x}_1}_{J_1},{\mathbf{x}_2}_{J_2})=(x_{j_{1,1}})=(x_2)$. Also,
\begin{equation}\nonumber
	\begin{split}
	f\arrowvert_{x_{2}=0}&=3x_1x_3+2x_4x_5+2,\\
	f\arrowvert_{x_{2}=1}&=3x_1x_3+2x_4x_5+2x_1+4x_3+x_4+x_5+3.
	\end{split}
\end{equation}
From (\ref{key114824}), we have
\begin{equation}\nonumber
	\begin{split}
	C_t=\left\{f+3(d_{1,1}+t_{1,1})x_2+3(d_{1,2}x_1+t_{1,2}x_3)+2(d_{2,1}x_4+t_{2,1}x_5):d_{1,1},~d_{1,2}\in \mbf Z_2,~d_{2,1}\in \mbf Z_3\right\}, 
	\end{split}
\end{equation}
where $t=4t_2+t_1$, $t_1=t_{1,1}+2t_{1,2}\in \mbf Z_2$ and $t_2=t_{2,1}\in \mbf Z_3$. It can be verified that the restricted functions $f\arrowvert_{x_{2}=0}$ and $f\arrowvert_{x_{2}=1}$ satisfy all the properties of \textbf{Corollary \ref{corr160824}}. Therefore, $\mal C=\{\psi(C_t):0\leq t\leq 11\}$ forms $(12,72)$-CCCs over the alphabet $\mal A_6=\{\zeta_6^i:i=0,1,\hdots,5\}$. 
Among $12$ codes, below we present two codes $\psi(C_1)$ and $\psi(C_{1,1})$, and their correlation properties. Similarly, the correlation properties of the other codes can be verified numerically. 
\begin{equation}\nonumber
	\begin{split}
\psi(C_1)=\begin{bmatrix}
\psi(f+3(0+1) x_2+3(0\cdot x_1+0\cdot x_3)+2(0\cdot x_4+0\cdot x_5))\\
\psi(f+3(0+1) x_2+3(0\cdot x_1+0\cdot x_3)+2(1\cdot x_4+0\cdot x_5))\\
\psi(f+3(0+1) x_2+3(0\cdot x_1+0\cdot x_3)+2(2\cdot x_4+0\cdot x_5))\\
\psi(f+3(1+1) x_2+3(0\cdot x_1+0\cdot x_3)+2(0\cdot x_4+0\cdot x_5))\\
\psi(f+3(1+1) x_2+3(0\cdot x_1+0\cdot x_3)+2(1\cdot x_4+0\cdot x_5))\\
\psi(f+3(1+1) x_2+3(0\cdot x_1+0\cdot x_3)+2(2\cdot x_4+0\cdot x_5))\\
\psi(f+3(0+1) x_2+3(1\cdot x_1+0\cdot x_3)+2(0\cdot x_4+0\cdot x_5))\\
\psi(f+3(0+1) x_2+3(1\cdot x_1+0\cdot x_3)+2(1\cdot x_4+0\cdot x_5))\\
\psi(f+3(0+1) x_2+3(1\cdot x_1+0\cdot x_3)+2(2\cdot x_4+0\cdot x_5))\\
\psi(f+3(1+1) x_2+3(1\cdot x_1+0\cdot x_3)+2(0\cdot x_4+0\cdot x_5))\\
\psi(f+3(1+1) x_2+3(1\cdot x_1+0\cdot x_3)+2(1\cdot x_4+0\cdot x_5))\\
\psi(f+3(1+1) x_2+3(1\cdot x_1+0\cdot x_3)+2(2\cdot x_4+0\cdot x_5))\\
\end{bmatrix},
	\end{split}
\end{equation}
and 
\begin{equation}
\begin{split}
\psi(C_{11})=\begin{bmatrix}
\psi(f+3(0+1) x_2+3(0\cdot x_1+1\cdot x_3)+2(0\cdot x_4+2\cdot x_5))\\
\psi(f+3(0+1) x_2+3(0\cdot x_1+1\cdot x_3)+2(1\cdot x_4+2\cdot x_5))\\
\psi(f+3(0+1) x_2+3(0\cdot x_1+1\cdot x_3)+2(2\cdot x_4+2\cdot x_5))\\
\psi(f+3(1+1) x_2+3(0\cdot x_1+1\cdot x_3)+2(0\cdot x_4+2\cdot x_5))\\
\psi(f+3(1+1) x_2+3(0\cdot x_1+1\cdot x_3)+2(1\cdot x_4+2\cdot x_5))\\
\psi(f+3(1+1) x_2+3(0\cdot x_1+1\cdot x_3)+2(2\cdot x_4+2\cdot x_5))\\
\psi(f+3(0+1) x_2+3(1\cdot x_1+1\cdot x_3)+2(0\cdot x_4+2\cdot x_5))\\
\psi(f+3(0+1) x_2+3(1\cdot x_1+1\cdot x_3)+2(1\cdot x_4+2\cdot x_5))\\
\psi(f+3(0+1) x_2+3(1\cdot x_1+1\cdot x_3)+2(2\cdot x_4+2\cdot x_5))\\
\psi(f+3(1+1) x_2+3(1\cdot x_1+1\cdot x_3)+2(0\cdot x_4+2\cdot x_5))\\
\psi(f+3(1+1) x_2+3(1\cdot x_1+1\cdot x_3)+2(1\cdot x_4+2\cdot x_5))\\
\psi(f+3(1+1) x_2+3(1\cdot x_1+1\cdot x_3)+2(2\cdot x_4+2\cdot x_5))\\
\end{bmatrix}.
\end{split}
\end{equation}
In Figure \ref{ankut111}-(a), we present the AACFs for the codes $\psi(C_{1})$ and $\psi(C_{11})$, and in \ref{ankut111}-(b), we present the ACCFs between them. 
\end{example}
Following the similar trend as \textbf{Theorem \ref{th314824}} and \textbf{Corollary \ref{corr160824}}, the proposed construction of $(p_1^{n_1+1}p_2^{n_2+1},L)$-CCCs can be extended to $(K,L)$-CCCs, where $K=\prod_{i=1}^k p_i^{n_i+1}$, $L=\prod_{i=1}^k p_i^{m_i}$ and $k\geq 3$. Again we recall the notations defined in Section \ref{seca16824} - \ref{sectn5}. For $\pi_1^c:\{1,2,\hdots,m_1-n_1\}\rightarrow \{1,2,\hdots,m_1\}\setminus J_1$ and $\pi_i^c:\{1,2,\hdots,m_i-n_i\}\rightarrow \{m_1+m_2+\cdots+m_{i-1}+1,\hdots,m_1+m_2+\cdots+m_i\}\setminus J_i$ are one-to-one mappings, where $\mbf c\in \mal V_{L'}$ is the vector representation of $c$, and $i=2,3,\hdots,k$. Again let $f:\mal V_L\rightarrow \mbb Z_q$ be a $q$-ary function such that after restricting the function $f$ at $\mbf x_J=\mbf c$, it reduces to a $q$-ary function of Hamming degree $2$ over the $m-n$ variables labelled $\{1,2,\hdots,m\}\setminus J$ as follows: 
\begin{equation}\label{func26824}
	\begin{split}
	f\arrowvert_{\mathbf{x}_J=\mathbf{c}}=&\sum_{i=1}^k\left(\frac{q}{p_i}\sum_{j=1}^{m_i-n_i-1} f_{i,j}(x_{\pi_{i}^c(j)})f_{i,j}'(x_{\pi_{i}^c(j+1)})+\sum_{j'=1}^{m_i-n_i} g_{i,j'}(x_{\pi_{i}^c(j')})\right) \\&+\sum_{i'=1}^{k-1} \lambda_{i'} f_{i'}(x_{\pi_{i'}^c(m_{i'}-n_{i'})}) h_{i'}(x_{\pi_{i'+1}^c(1)}),
	\end{split}
\end{equation}
where $f_{i,j},~f_{i,j}',~g_{i,j'},~f_{i'}$, and $h_{i'}$ are functions from $\mbb Z_q$ to $\mbb Z_q$, and $\lambda_{i'}\in \mbb Z_q$. Let us again define $t=t_1+\sum_{\beta=2}^k t_\beta \prod_{\alpha=1}^{\beta-1} p_{\alpha}^{n_\alpha+1}$, where $t_u=\sum_{v=1}^{n_u+1}t_{u,v}p_u^{v-1}$ and $u=1,2,\hdots,k$. Then for $0\leq t<K$, we define $C_t$ as the set of $K$ functions as follows:
\begin{equation}
	\begin{split}
	C_t=&\left\{ f+ \sum_{i=1}^k\frac{q}{p_i}\left(\left(\mbf d_i+\mbf t_i\right)\cdot {\mbf x_i }_{J_i}+\left(d_{i,n_i+1}x_{\pi_i(1)}+t_{i,n_i+1}x_{\pi_i(m_i-n_i)} \right)\right):\right.\\&~~~~ \left.\mbf d_i\in \mbf Z_{p_i}^{n_i}, d_{i,n_i+1}\in \mbf Z_{p_i}, i=1,2,\hdots,k  \right\},
	\end{split}
\end{equation}
where $\mbf t_i=(t_{i,1},t_{i,2},\hdots,t_{i,n_i})$ and $\mbf d_i=(d_{i,1},d_{i,2},\hdots,d_{i,n_i})$.
\begin{corollary}\label{corr19824}
For any choices of $g_{i,j'},~f_{i'},~h_{i'}$, one-to-one mapping $\pi_i^c$, and for any values of $\lambda_{i'}\in\mbb Z_q$ in (\ref{func26824}), the set of codes $\mal C=\{\psi(C_t):0\leq t<K\}$ forms $(K,L)$-CCCs iff 
$f_{i,j}$ and $f_{i,j}'$ permute the set $\mbf Z_{p_i}$ under modulo $p_i$, where $i=1,2,\hdots,k$, $j=1,2,\hdots,m_i-n_i-1$, $j'=1,2,\hdots,m_i-n_i$, and $i'=1,2,\hdots,k-1$. 
\end{corollary}
\textbf{Corollary \ref{corr19824}} can be proved using a similar technique as introduced in \textbf{Theorem \ref{th314824}}.
\begin{remark} 
In (\ref{func26824}), by setting \(\lambda_i' = 0\), the function \(f\arrowvert_{\mathbf{x}_J=\mathbf{c}}\) can be represented as the direct sum of \(k\) functions, denoted as \(F_1, F_2, \dots, F_k\), where 
\[
F_i = \frac{q}{p_i} \sum_{j=1}^{m_i-n_i-1} f_{i,j}(x_{\pi_{i}^c(j)})f_{i,j}'(x_{\pi_{i}^c(j+1)}) + \sum_{j'=1}^{m_i-n_i} g_{i,j'}(x_{\pi_{i}^c(j')}).
\]
By setting \(q = p_i\) and \(f = F_i\) in \textbf{Corollary \ref{cor15824}}, one can derive \((p_i^{n_i+1}, p_i^{m_i})\)-CCCs over the alphabet \(\mathcal{A}_{p_i}\) iff \(f_{i,j}\) and \(f_{i,j}'\) permute \(\mathbb{Z}_{p_1}\). Then, by taking the Kronecker product of these \(k\) CCCs, i.e., \((p_i^{n_i+1}, p_i^{m_i})\)-CCCs for \(i = 1, 2, \dots, k\), one can construct \((K, L)\)-CCCs over the alphabet \(\mathcal{A}_q\). Since \(f\arrowvert_{\mathbf{x}_J=\mathbf{c}} = \sum_{i=1}^k F_i\), it follows that 
\[
\psi(f\arrowvert_{\mathbf{x}_J=\mathbf{c}}) = \psi(F_1) \otimes \psi(F_2) \otimes \cdots \otimes \psi(F_k).
\]
Thus, by setting \(\lambda_i = 0\) in \textbf{Corollary \ref{corr19824}}, one can also obtain \((K, L)\)-CCCs expressed as the Kronecker product of \(k\) CCCs defined by \(F_1, F_2, \dots, F_k\). This shows that using \textbf{Corollary \ref{corr19824}}, one can directly obtain functional representations for these \((K, L)\)-CCCs, which can be derived by performing the Kronecker product on \((p_i^{n_i+1}, p_i^{m_i})\)-CCCs constructed from \textbf{Corollary \ref{cor15824}}.
\end{remark} 

\begin{remark}[Comparison with \cite{sarkar2021multivariable}]
By choosing $f_{i,j}(x_{\pi_{i}^c(j)})=x_{\pi_{i}^c(j)}$, $f_{i,j}'(x_{\pi_{i}^c(j+1)})=x_{\pi_{i}^c(j+1)}$, $\lambda_i'=0$, and $g_{i,j'}$ as a liner function, from 
\textbf{Corollary \ref{corr19824}}, we obtain $(K,L)$-CCCs which is covered in \cite{sarkar2021multivariable}. Hence, the result in \cite{sarkar2021multivariable} appears as a special case of our proposed work. 	
\end{remark}
%%%%%%%%%%%%%%%%%%%%%%%%%%%%%%%%%%%%%%%%%%%%%%%%%%%%%%%%%%%%%%%%%%%%%
%%%%%%%%%%%%%%%%%%%PREVIOUS CONTRIBUTION ENDS HERE%%%%%%%%%%%%%%%

In the following remarks, we compare our proposed results with existing related results.
% intorduced in \cite{wang2020new,wangong}, and with \cite{shen2023}.
	\begin{remark}[Comparison with \cite{Davis1999, pater2000, rati}]
		In \textbf{Corollary \ref{corr19824}}, the variable $K$ signifies the number of CSs within the set of CCCs, including their constituent sequences, while $L$ represents the sequence length. Additionally, $\mathbb{Z}_q$ denotes the alphabet set. We can revisit the expressions for $K$, $L$, and $q$ as follows:
		\[
		K=\prod_{i=1}^k p_i^{n_i+1},~L=\prod_{i=1}^kp_i^{m_i},~q=\prod_{i=1}^k p_i.
		\]
		%where $n$ can be chosen as any positive integer. However, opting for $n=1$ results in the proposed framework generating sequences with the smallest %alphabet size.
		 For $k=1$ and $p_1=2$, the parameters can be expressed as:
		\[
		K=2^{n_1+1},~L=2^{m_1},~q=2,
		\] 
		underscoring that the findings reported in \cite{Davis1999, pater2000, rati} with the aforementioned parameters emerge as special cases of \textbf{Corollary \ref{corr19824}}, where \cite{Davis1999, pater2000, rati} introduce the construction of GCPs, CSs, and CCCs with the parameters as mentioned above. It is to be noted here that our proposed construction also works if we choose $q$ to be a positive integer such that $p_i|q$, $i=1,2,\hdots,k$. As we are interested in a small alphabet, throughout the paper we fix $q$ to its minimum value which is $\prod_{i=1}^k p_i$. 
\end{remark} 
\begin{remark}[Comparison with \cite{wang2020new,wangong}]
In \cite{wang2020new}, PU matrix-based construction of CCCs has been reported, followed by a systematic approach to
extract $p$-ary functions, for prime $p$, as the generators for CCCs of non-power-of-two lengths.
%A multivariable function can generate CCCs directly as it does not require to perform a series of operations on
%any special tools like the proposed construction.
The construction method in \cite{wang2020new} needs to perform the following steps {(unlike our proposed one)}
to generate CCCs:
\begin{itemize}
\item It needs to follow a generic framework for constructing desired PU matrices.
\item A series of operations need to be performed to extract $p$-ary functions from the constructed PU matrices.
\end{itemize}
In \cite{wangong}, Wang-Gong proposed a construction of $(N,N^m)$-CCCs using $N$-ary functions, where $N$ is an integer power, which are obtained through the
	following steps:
	\begin{itemize}
	\item First, construct a Butson-type Hadamard (BH) matrix of order $N$ with the help of $N$-ary sequences with two-level autocorrelation of period $N-1$.
	\item Then, construct $\delta$-quadratic terms with the help of the constructed BH matrix of order $N$.
	\item Finally, with the help of $\delta$-quadratic terms and suitable permutation polynomials, $N$-ary functions are obtained to
	generate $(N,N^m)$-CCCs over the alphabet $\mathbb{Z}_{N}$.
	\end{itemize}
It can be easily verified that the proposed constructions in this paper generate more flexible parameters as compared to \cite{wang2020new} and  \cite{wangong}.
\end{remark}
\begin{remark}[Comparison with \cite{shen2023}]
In the context of the current constructions, the importance of defining the domain and co-domain of a function becomes apparent when considering the alphabet size of the resultant sequences. This study characterizes the domain as $\calV_L \subset \mathbb{Z}_q^m$ and the co-domain as $\mathbb{Z}_q$. This configuration generates sequences of length $L$ over the alphabet $\mathbb{Z}_q$, which is determined by the prime factors of $L$. Notably, when co-primes values for $m_i \neq m_j$ occur for some specific indices $i \neq j \in \{1, 2, \ldots, k\}$, the proposed construction leads to the creation of CCCs featuring the minimal alphabet $q=\prod_{i=1}^k p_i$. This outcome deviates from existing constructions, particularly in cases where $k \geq 2$. For instance, in a recent work \cite{shen2023}, the authors adopt a domain $\mathbb{Z}_q^m$, where $m \geq 2$, and a co-domain of $\mathbb{Z}_q$, yielding $(q^{n+1}, q^m)$-CCCs with an alphabet size of $q$, where $0\leq n\leq m-1$. When dealing with scenarios where $L = \prod_{i=1}^k p_i^{m_i}$ and some of  $m_1,m_2,\hdots,m_k$ are co-primes to each-other, \cite{shen2023} can also generate a sequence length of $L$ by setting $m = 1$ and $q = L$. However, this yields both the alphabet size and set size as $L$, which is relatively larger than our outcomes. In practical applications, CCCs characterized by modest set and alphabet sizes tend to outperform those with larger dimensions, owing to considerations like the PMEPR of constituent sequences \cite{pater2000, Davis1999}. This is where our proposed construction showcases its strengths, achieving superior performance when scrutinized through vital parameters such as alphabet and set size. 
\end{remark}
%%%%%%%%%%%%%Conclusion%%%%%%%%%%%%%%%%
%%%%%%%%%%%%%%%%%%%%%%%%%%%%%%%%%%%%%%%
\section{Conclusion}\label{conclu200924}
In connection with $q$-ary CCCs, unlike the existing works most of which contribute to only sufficient conditions that provide some specific forms for the functional representations of CCCs, in this paper, we have studied both the necessary and sufficient conditions on the $q$-ary functional representations for $q$-ary CCCs. This study on necessary conditions contributes in understanding the complete relationship between $q$-ary CCCs and its functional representations in the form as presented in (\ref{func26824}) and (\ref{qry_n1}). Please note that in (\ref{func26824}) and (\ref{qry_n1}), the properties of desired $q$-ary functions has been presented in terms of its restrictions over some variable, and through all possible restrictions of a function, we can identify the main function easily. First, We have studied the $q$-ary functions on the domain $\mbb Z_q^m$ where $q\geq 2$ is any integer, and determine the complete class of functions for $(q^{n+1},q^m)$-CCCs over an alphabet of size $q$. Though the $q$-ary functions on the domain $\mbb Z_q^m$, one can produce almost all possible lengths CCCs, but in many cases such as CCCs of lengths in the form $\prod_{i=1}^k p_i^{m_i}$, where $p_1,p_2,\hdots,p_k$ are  distinct primes, $m=m_1+m_2+\hdots+m_k$, and $m_i$ is co-prime to $m_j$ for some $1\leq i,j\leq k$, the set size and alphabet size become as large as sequence lengths. To overcome this issue, we propose another type of $q$-ary functions where $q=\prod_{i=1}^kp_i$, and the domain appears in the form $\mbf Z_{p_1}^{m_1}\times\mbf Z_{p_2}^{m_2}\times \cdots\times \mbf Z_{p_k}^{m_k}\subset \mbb Z_q^m$. This type of $q$-ary functions generate $(\prod_{i=1}^k p_i^{n_i+1},\prod_{i=1}^k p_i^{m_i})$-CCCs over an alphabet of size $\prod_{i=1}^kp_i$. 
%%%%%%%%%%%%%%%%%%%%%%%%%%%%%%%%%%%%%%%%%%%%%%%%%%%%%%%%%%%%%%%%%%%%%%%%%%%%%%%%
\appendices
\section{Details of the Correspondence Between $\mbf x$ and $x$}\label{appendix1}
We define necessary algebraic operations on the set $\mal{V}_L$  in order to introduce $q$-ary functions while maintaining $\mal{V}_L$ as the domain. 
We start with a given positive integer $L = L_1L_2$, define a mapping $\rho$: $\mathbf{Z}_L \rightarrow \mathbf{Z}_{L_1}\times \mathbf{Z}_{L_2}$ as 
$\rho(x):= (\sigma_1, \sigma_2)$ where 
\[
\sigma_1 = x \bmod{L_1} \text{ and } \sigma_2 = (x-\sigma_1)/L_1.
\] It is clear that $\rho$ is bijective and its inverse $\rho^{-1}$ is given by \[x = \rho^{-1}(\sigma_1, \sigma_2) = \sigma_2*L_1 + \sigma_1.\]
More generally, for an integer $L=\prod_{i=1}^k L_i$,  we can define a one-to-one mapping 
$\rho$: $\mathbf{Z}_L \rightarrow \mathbf{Z}_{L_1}\times \mathbf{Z}_{L_2}\times \dots \times \mathbf{Z}_{L_k}$ as $\rho(x) = (\sigma_1, \sigma_2, \hdots, \sigma_k)$, where $\sigma_i$'s is recursively derived as follows: 
\begin{equation}\label{Eq_Int2Vec}
\begin{split}
& y_1 = x \\
& \sigma_t = y_t \bmod{L_t} \text{ and } y_{t}= \frac{y_{t-1}-\sigma_{t-1}}{L_{t-1}} \quad \text{ for } t =2,3,\hdots,k. 
\end{split}
\end{equation}
From the above recursive relations, it is easy to verify that the pre-image of $(\sigma_1,\sigma_2,\dots, \sigma_k)$ is given by 
\begin{equation}\label{Eq_Vec2Int}
\begin{split}
x=\rho^{-1}((\sigma_1,\sigma_2, \hdots, \sigma_k)) &= (y_{2}L_1 + \sigma_1) 
\\&= ((y_{3}L_{2} + \sigma_{2})L_1 + \sigma_1) =\cdots =
%\\&= ((\sigma_1L_{2}+\sigma_2)L_3 + \sigma_3)L_4\cdots)L_k + \sigma_k
\\&=\sigma_k L_1 L_2\cdots L_{k-1} + \sigma_{k-1} L_1L_2\cdots L_{k-2} + \dots + \sigma_{2}L_{1} + \sigma_1
\\&= \sigma_1 \Delta_1 + \sigma_2 \Delta_2 + \dots + \sigma_k \Delta_k,
\end{split}
\end{equation} where $\Delta_i = \frac{L}{L_iL_{i+1}\cdots L_{k}}$. Note that the above relation does not pose any restriction on the relation among $L_1, L_2,\dots, L_k$. In particular, when $L_1=L_2=\dots = L_k = p$, the image $\rho(x)=(\sigma_1,\sigma_2,\hdots, \sigma_k)$ corresponds to the conventional $p$-ary expansion of $x = \sum_{i=1}^{k}\sigma_i\Delta_i$, where $\Delta_i = \frac{p^k}{\prod_{i\leq j\leq k}p} = p^{i-1}$. 
%Suppose a positive integer $L$ has a decomposition $L = L_1L_2\dots L_k$ where $L_i = p_i^{m_i}$ and $p_1<p_2<\dots <p_k$.
From this discussion, we can define an one-to-one mapping $\rho: \mathbf{Z}_L \rightarrow \mal{V}_L$ as 
$$\mathbf{x} := \rho(x) = (\bx_1, , \dots, \bx_k) = (x_1, \dots, x_{m_1}, \dots, x_{m_1+\cdots+m_{k-1}+1}, \dots, x_{m}),
$$ where 
\begin{equation}\nonumber
\begin{split}
\mbf x_1&=(x_1,x_2,\hdots,x_{m_1})~\in \mbf Z_{p_1}^{m_1},~\textnormal{and}\\
\bx_i& = (x_{m_1+\cdots+m_{i-1}+1}, \dots, x_{m_1+\hdots+m_i}) \in \bZ_{p_i}^{m_i},~2\leq i\leq k.
\end{split}
\end{equation}
%$\mbf x_1=(x_1,x_2,\hdots,x_{m_1})$, and $\bx_i = (x_{m_1+\cdots+m_{i-1}+1}, \dots, x_{m_1+\hdots+m_i}) \in \bZ_{p_i}^{m_i}$, $2\leq i\leq k$. 
Also, 
\begin{equation}\label{ref100}
\begin{split}
\sigma_1&=\sum_{j=1}^{m_1}x_jp_1^{j-1},~\textnormal{and}\\
\sigma_i&=\sum_{j=1}^{m_i} x_{m_1+\cdots+m_{i-1}+j} p_i^{j-1},~i=2,3,\hdots,k.
\end{split}
\end{equation}
% is the $p_i$-ary expansion of the corresponding integer $x_i$ in $\bZ_{p_i^{m_i}}$, namely, $x_i = \sum\limits_{j=1}^{m_i}x_{i,j}p_i^{m_i-j}$. 
Let $\bx$ is the vectorial representation of an integer $x \in \bZ_L$ and $x$ is the integer representation of a vector $\bx \in \mal{V}_L$, which are connected via the following equality
$$x=\rho^{-1}(\mbf x)=\sigma_1\Delta_1+\sigma_2\Delta_2+\hdots+\sigma_k\Delta_k,$$
where $\sigma_i$ is presented in (\ref{ref100}) for $i=1,2,\hdots,k$.
\bibliographystyle{IEEEtran}
\bibliography{CCC_SNC_CCC}
\end{document}